\documentclass[10pt]{article}
\usepackage{amssymb,amsmath,amsthm}
\usepackage{graphicx}

 \UseRawInputEncoding
 
\begin{document}
 \title{\bf Ground state solutions for a non-local type problem in fractional Orlicz Sobolev
spaces}
 \author { Liben Wang$^{1}$\footnote{Corresponding author, E-mail address: wanglbdgust@163.com.},\  \ Xingyong Zhang$^{2}$ \ and Cuiling Liu$^2$\\
       {\footnotesize   $^1$School of Computer Science and Technology, Dongguan University of Technology,}\\
        {\footnotesize  Dongguan, Guangdong, 523808, P.R. China}\\
       {\footnotesize   $^2$Faculty of Science, Kunming University of Science and Technology,}\\
        {\footnotesize  Kunming, Yunnan, 650500, P.R. China}\\}
 \date{}
 \maketitle{}
 \begin{center}
 \begin{minipage}{15cm}
 \par
 \small  {\bf Abstract:}  In this paper, we study the following nonlocal problem in fractional Orlicz Sobolev spaces
 \begin{eqnarray*}
 (-\Delta_{\Phi})^{s}u+V(x)a(|u|)u=f(x,u),\quad x\in\mathbb{R}^N,
 \end{eqnarray*}
 where $(-\Delta_{\Phi})^{s}(s\in(0, 1))$ denotes the non-local and maybe non-homogeneous operator, the so-called fractional $\Phi$-Laplacian. Without assuming the Ambrosetti-Rabinowitz type and the Nehari type conditions on the nonlinearity, we obtain the existence of ground state solutions for the above problem in periodic
case. The proof is based on a variant version of the mountain pass theorem and a Lions' type result for fractional Orlicz Sobolev spaces.
 \par
 {\bf Keywords:} Fractional Orlicz-Sobolev spaces; Fractional $\Phi$-Laplacian; Critical point; Ground state\\
 {{ \bf 2010 Mathematics Subject Classification:} 35R11, 46E30, 35A15}

 \end{minipage}
 \end{center}
 \section{Introduction and main results}
 \setcounter{equation}{0}
 \allowdisplaybreaks

  \ \quad In the last decades, much attention has been devoted to the study of the nonlinear Schr\"{o}dinger equations involving non-local operators. These types of operators can be used to model many phenomena in the natural sciences such as fluid dynamics, quantum mechanics, phase transitions, finance and so on, see \cite{Laskin2000,Alberti1998,Metzler2004,Mosconi2016} and the references therein. Due to the important work of Fern\'andez Bonder and Salort \cite{Bonder2019}, a new generalized fractional $\Phi$-Laplacian operator has caused great interest among scholars in recent years since it allows to model non-local problems involving a non-power behavior, see \cite{Bonder2023,Salort2022,Salort2020,Alberico2021,Azroul2020,Bahrouni2020,Chaker2022,Silva2021} and the references therein.

  In this paper, we are interested in studying the following nonlocal problem involving fractional $\Phi$-Laplacian:
 \begin{eqnarray}\label{eq1}
 (-\Delta_{\Phi})^{s}u+V(x)a(|u|)u=f(x,u),\quad x\in\mathbb{R}^N,
 \end{eqnarray}
 where $s\in(0, 1), N\in\mathbb{N}$, the function $a:[0,+\infty)\rightarrow\mathbb{R}$ is such that
 $\phi:\mathbb{R}\rightarrow\mathbb{R}$ defined by
\begin{eqnarray}\label{1.2}
 \phi(t)=
 \begin{cases}
  \begin{array}{ll}
 a(|t|)t &\mbox{for } t\neq 0,\\
 0 &\mbox{for } t=0,\\
    \end{array}
 \end{cases}
 \end{eqnarray}
 is an increasing homeomorphism from $\mathbb{R}$ onto $\mathbb{R}$, $\Phi:[0,+\infty)\rightarrow[0,+\infty)$ defined by
 $$\Phi(t)=\int_{0}^{t}\phi(\tau)d\tau$$
 is an $N$-function (see Section 2 for details), which together with the potential $V$ and the nonlinearity $f$ satisfy the following basic assumptions:\\
 $(\phi_1)$\
 $1<l:=\inf_{t>0}\frac{t\phi(t)}{\Phi(t)}\leq\sup_{t>0}\frac{t\phi(t)}{\Phi(t)}=:m<\min\{\frac{N}{s},l^*\},$
 where $l^*:=\frac{Nl}{N-sl},$\\
 $(V)$  $V\in C(\mathbb{R}^N, \mathbb{R}_+)$ is 1-periodic in $x_1, \cdots, x_N$ (called 1-periodic in $x$ for short) and so there exist two constants $\alpha_1, \alpha_2>0$ such that $\alpha_1\leq V(x)\leq \alpha_2 \mbox{ for all }x\in \mathbb{R}^N;$\\
 $(f_1)$ $f\in C(\mathbb{R}^N\times\mathbb{R})$ is 1-periodic in $x$ satisfying
 $$\lim_{|t|\rightarrow 0}\frac{f(x,t)}{\phi(|t|)}=0\quad \mbox{ and } \quad
 \lim_{|t|\rightarrow \infty}\frac{f(x,t)}{\Phi_{*}'(|t|)}=0, \quad \mbox{ uniformly in }x\in \mathbb{R}^N,$$
 where $\Phi_{*}$ denotes the Sobolev conjugate function of $\Phi$ (see Section 2 for details).\\
 For $s\in(0, 1)$, the so-called fractional $\Phi$-Laplacian operator is defined as
 \begin{eqnarray}\label{1.3}
 (-\Delta_{\Phi})^{s}u(x):=P.V.\int_{\mathbb{R}^N}a(|D_s u|)\frac{D_s u}{|x-y|^{N+s}}dy,\quad \mbox{ where } \quad D_s u:=\frac{u(x)-u(y)}{|x-y|^{s}}
 \end{eqnarray}
 and $P.V.$ denotes the principal value of the integral. Notice that if $\Phi(t)=|t|^p (p>1)$, then fractional $\Phi$-Laplacian operator reduces to the following fractional $p$-Laplacian operator
 $$(-\Delta_{p})^{s}u(x):=P.V.\int_{\mathbb{R}^N}\frac{|u(x)-u(y)|^{p-2}(u(x)-u(y))}{|x-y|^{N+ps}}dy.$$
 To study this class of nonlocal problem involving fractional $p$-Laplacian, variational method has become one of the important tools over the past several decades, see \cite{Dipierro2013,Chang2013,Secchi2016,Nezza2012,Ambrosio2018,Perera2016,Xu2015} and the references therein. In most of the references, to ensure the boundedness of the Palais-Smale sequence or Cerami sequence of the energy functional, the following (AR) type condition for the nonlinearity $f$ due to Ambrosetti-Rabinowitz \cite{Ambrosetti1973} was always assumed:\\
 (AR) there exists $\mu>p$ such that
 $$0<\mu F(x,t)\leq tf(x,t), \quad\mbox{ for all } t\neq0,$$
 where and in the sequel, $F(x,t)=\int_{0}^{t}f(x,\tau)d\tau$. \\
 (AR) implies that there exist two positive constants $c_1, c_2$ such that
 $$F(x,t)\geq c_1|t|^{\mu}-c_2, \quad \mbox{ for all } (x,t)\in\mathbb{R}^N\times\mathbb{R},$$ which is obviously stronger than the following $p$-superlinear growth condition:\\
 $(F_1)$ $\lim_{|t|\rightarrow\infty}\frac{F(x,t)}{|t|^p}=+\infty$, uniformly in $x\in\mathbb{R}^N.$\\
$(F_1)$ was first introduced by Liu and Wang in \cite{Liu2004} for the case $p=2$ and then was commonly used in recent papers. With the development of variational theory and application, some new restricted conditions were established in order to weaken (AR). But, most of those conditions are just complementary to (AR). For example, replace (AR) with $(F_1)$ and the following Nehari type condition:\\
(Ne) $\frac{f(x,t)}{|t|^{p-1}}$ is (strictly) increasing in $t$ for all $x\in\mathbb{R}^N$.\\
For the case $p=2$, Li, Wang and Zeng proved the existence of ground state by Nehari method in \cite{Li2006}.
Besides, for the case $p=2$, Ding and Szulkin in \cite{Ding2007} replaced (AR) with $(F_1)$ and the following condition:\\
$(F_2)$ $\mathcal{F}(x,t)>0$ for all $t\neq 0$, and $|f(x,t)|^{\sigma}\leq c_3\mathcal{F}(x,t)|t|^{\sigma}$ for some $c_3>0, \sigma>\max\{1,\frac{N}{2}\}$ and all $(x,t)$ with $|t|$ larger enough, where $\mathcal{F}(x,t)=tf(x,t)-2F(x,t)$.\\
 They proved that $(F_1)$ and $(F_2)$ hold if the nonlinearity $f$ satisfies (AR) and a subcritical growth condition that $|f(x,t)|\leq c_4(|t|+|t|^{q-1})$ for some $c_4>0, q\in(2, 2^*)$ and all $(x,t)\in\mathbb{R}^N\times\mathbb{R}$, where $2^*=\frac{2N}{N-2}$ if $N\geq3$ and $2^*=\infty$ if $N=1$ or $N=2$. Some conditions similar to $(F_2)$ were introduced in \cite{Lin2013,Cheng2016} for the case $p>1$. Moreover, in \cite{Tang2014}, Tang introduced the following new and weaker super-quadratic condition:\\
 $(F_3)$ there exists a $\theta_0\in(0, 1)$ such that
 $$\frac{1-\theta^2}{2}tf(x,t)\geq\int_{\theta t}^{t}f(x,\tau)d\tau=F(x,t)-F(x,\theta t),\mbox{ for all } \theta\in[0, \theta_0], (x,t)\in\mathbb{R}^N\times\mathbb{R}.$$
 Tang proved that $(F_3)$ is weaker than both (AR) and (Ne) and also different from $(F_2)$. It's worth noting that $(F_3)$ has been extended for the case $p>1$ in \cite{Mi2021}.
 \par
 To the best of our knowledge, some conditions mentioned above have been successfully generalized to the nonlocal problem involving fractional $\phi$-Laplacian. In \cite{Sabri2023}, for equation (\ref{eq1}) with potential $V(x)\equiv1$, by applying the mountain pass theorem, Sabri-Ounaies-Elfalah proved the existence of nontrivial solution when the autonomous nonlinearity $f(u)$ satisfies an (AR) type condition. On the whole space $\mathbb{R}^N$, to overcome the difficulty due to the lack of compactness of the Sobolev embedding, the authors reconstructed the compactness by choosing a radially symmetric function subspace as the working space. In \cite{Silva2021}, for equation (\ref{eq1}) with unbounded or bounded potentials $V$, by applying the Nehari manifold method, Silva-Carvalho-de Albuquerque-Bahrouni proved the existence of ground state solutions when the nonlinearity $f$ satisfies the following both (AR) and (Ne) type conditions:\\
 $\mbox{(AR)}^*$ there exists $\theta>m$ such that $\theta F(x,t)\leq tf(x,t)$, for $(x,t)\in\mathbb{R}^N\times\mathbb{R};$\\
 $\mbox{(Ne)}^*$  the map $t\rightarrow\frac{f(x,t)}{|t|^{m-1}}$ is strictly increasing for $t>0$ and strictly decreasing for $t<0$.\\
 To be precise, for the case when $V$ is unbounded, the authors reconstructed the compactness by assuming that $V$ is coercive and then choosing a subspace depending on $V$ as the working space. For the case when $V$ is bounded, to overcome the difficulty due to the lack of compactness and obtain a nontrivial solution, the authors assumed that $V$ and $f$ are 1-periodic in $x$ and introduced an important Lions' type result for fractional Orlicz-Sobolev spaces (see Theorem 1.6 in \cite{Silva2021}). Since the ground state solution there is obtained as a minimizer of the energy functional on the Nehari manifold $\mathcal{N}$, it is crucial to require that $f$ is of class $C^1$. Otherwise $\mathcal{N}$ may not be a $C^1$-manifold and it is not clear that the minimizer on the Nehari manifold $\mathcal{N}$ is a critical point of the energy functional.
 \par
Motivated by \cite{Silva2021}, in this paper we still study the existence of ground state for equation (\ref{eq1}) under the assumption that $V$ and $f$ are 1-periodic in $x$. We manage to extend the above $p$-superlinear growth conditions $(F_2)$ and $(F_3)$ to the nonlocal problem involving fractional $\Phi$-Laplacian. Instead of applying the Nehari manifold method, we firstly prove that equation (\ref{eq1}) has a nontrivial solution by using a variant mountain pass theorem (see Theorem 3 in \cite{Silva2010}). Subsequently, we prove the existence of ground state by using the Lions' type result for fractional Orlicz-Sobolev spaces and some techniques of Jeanjean and Tanaka (see Theorem 4.5 in \cite{Jeanjean2002}).
 \par
 Next, we present our main results as follows.
\vskip2mm
 \noindent
 {\bf Theorem 1.1.}\ \  {\it  Assume that $(\phi_1)$, $(V)$, $(f_1)$ and the following conditions hold:\\
 $(\phi_2)$ $\limsup_{t\rightarrow 0}\frac{|t|^{l}}{\Phi(|t|)}<+\infty;$\\
 $(f_2)$ $\lim_{|t|\rightarrow \infty}\frac{F(x,t)}{\Phi(|t|)}=+\infty,$ uniformly in $x\in \mathbb{R}^N;$\\
 $(f_3)$ $\widehat{F}(x,t)>0$ for all $t\neq 0$, and $|F(x,t)|^{k}\leq c\widehat{F}(x,t)|t|^{lk}$ for some $c>0,$ $k>\frac{N}{sl}$ and all $(x,t)$ with $|t|$ larger enough, where $\widehat{F}(x,t)=tf(x,t)-mF(x,t).$\\
 Then equation (\ref{eq1}) has at least one ground state solution. }
 \vskip2mm
 \noindent
 {\bf Theorem 1.2.}\ \  {\it  Assume that $(\phi_1)$, $(V)$, $(f_1)$ and the following conditions hold:\\
  $(f_4)$ $F(x,t)\geq0$ for all $(x,t)\in\mathbb{R}^N\times\mathbb{R},$ and $\lim_{|t|\rightarrow \infty}\frac{F(x,t)}{|t|^m}=+\infty,$ uniformly in $x\in \mathbb{R}^N;$\\
 $(f_5)$ there exists a $\theta_0\in(0, 1)$ such that
 $$\frac{1-\theta^l}{m}tf(x,t)\geq\int_{\theta t}^{t}f(x,\tau)d\tau=F(x,t)-F(x,\theta t), \mbox{ for all } \theta\in[0, \theta_0], (x,t)\in\mathbb{R}^N\times\mathbb{R}.$$
 Then equation (\ref{eq1}) has at least one ground state solution. }
 \vskip2mm
\noindent
 {\bf Remark 1.3.}\ \   To some extent, Theorem 1.2 improves the result of Theorem 1.8 in \cite{Silva2021}. In fact, our results do not require the smoothness condition that functions $f$ and $a$ are of class $C^1$. Moreover, it is obvious that $(\varphi_4)$ in \cite{Silva2021} implies $(\phi_1)$ and $(f_0)$ in \cite{Silva2021} implies our subcritical growth condition given by $(f_1)$. Furthermore, when $\Phi(t)=|t|^2$, $(f_5)$ is weaker than both (AR) type condition $(f_4)$ and (Ne) type condition $(f_4)$ in \cite{Silva2021} (see \cite{Tang2014}).
 \vskip2mm
\noindent
 {\bf Remark 1.4.}\ \  Theorem 1.2 extends and improves the result of Theorem 1.1 in \cite{Zhang2017}. In fact, when $\Phi(t)=|t|^2$, our subcritical growth condition given by $(f_1)$ reduces to
 \begin{eqnarray}\label{1.4}
 \lim_{|t|\rightarrow \infty}\frac{f(x,t)}{|t|^{2^*-1}}=0, \quad \mbox{ uniformly in }x\in \mathbb{R}^N,
 \end{eqnarray}
 which is weaker than $(A_2)$ in \cite{{Zhang2017}}. For example, it is easy to check that function $f(t)=\frac{|t|^{2^*-2}t}{\log(e+|t|)}$ satisfies (\ref{1.4}) but does not satisfy $(A_2)$ in \cite{{Zhang2017}}. Moreover, it is obvious that Theorem 1.1 is different from Theorem 1.2 in \cite{{Zhang2017}} even when the fractional $\Phi$-Laplacian equation (\ref{eq1}) reduces to the fractional Schr\"{o}dinger equation.
\par
The rest of this paper is organized as follows. In Section 2, we recall some definitions and basic properties on the Orlicz and fractional Orlicz Sobolev spaces. In Section 3, we complete the proofs of the main results. In Section 4, we present some examples about the function $\phi$ defined by (\ref{1.2}) and nonlinearity $f$ to illustrate our results.

 \vskip-2mm
 \noindent
 \section{Preliminaries}
 \setcounter{equation}{0}
 \vskip0mm
  \qquad  In this section, we make a brief introduction about Orlicz and fractional Orlicz Sobolev spaces. For more details, we refer the reader to \cite{Bonder2019,Adams2003,M. M. Rao2002} and references therein.
  \par
   First of all, we recall the notion of $N$-function. Let $\phi: [0,+\infty)\rightarrow [0,+\infty)$ be a right continuous, monotone increasing function with \\
 $1)$ $\phi(0)=0$;\\
 $2)$ $\lim\limits_{t\rightarrow +\infty} \phi (t)=+\infty$;\\
 $3)$ $\phi(t)>0 $ whenever $t>0$.\\
 Then the function defined on $[0, +\infty)$ by $\Phi(t)=\int_{0}^{t}\phi(\tau)d\tau$
 is called an $N$-function. It is obvious that $\Phi(0)=0$ and $\Phi$ is strictly increasing and convex in $[0, +\infty)$.
 \par
 An $N$-function $\Phi$ satisfies a $\Delta_2$-condition if there exists a constant $K>0$ such that $\Phi(2t)\leq K\Phi(t)$ for all $t\geq0$. $\Phi$ satisfies a $\Delta_2$-condition if and only if for any given $c\geq 1$, there exists a constant $K_c>0$ such that $\Phi(ct)\leq K_c\Phi(t)$ for all $t\geq0$.
 \par
Given two $N$-functions $A$ and $B$, we say that $B$ is essentially stronger than $A$ (or equivalently
that A decreases essentially more rapidly than B), denoted by $A\prec\prec B$, if for each $c>0$ there holds that
$$\lim_{t\rightarrow +\infty}\frac{A(ct)}{B(t)}=0.$$
  \par
    For $N$-function $\Phi$, the complement of $\Phi$ is given by
 $$\widetilde{\Phi}(t)=\max_{\rho\geq 0}\{t\rho-\Phi(\rho)\},\quad\mbox{ for } t\geq 0.$$
 Then, we have the Young's inequality, that is
 \begin{eqnarray}\label{ 2.11}
 \rho t\leq \Phi (\rho)+\widetilde{\Phi}(t), \quad\mbox{ for all } \rho, t\geq 0,
 \end{eqnarray}
 and the following inequality (see Lemma A.2 in \cite{Fukagai2006}), that is
 \begin{eqnarray}\label{ 2.1}
 \widetilde{\Phi}(\phi(t))\leq \Phi(2t), \quad\mbox{ for all } t\geq 0.
 \end{eqnarray}
\par
Now we recall the Orlicz space $L^{\Phi}(\mathbb{R}^N)$ associated with $\Phi$. When $\Phi$ satisfies $\Delta_2$-condition, the Orlicz space $L^{\Phi}(\mathbb{R}^N)$ is the vectorial space of the measurable functions $u: \mathbb{R}^N\rightarrow \mathbb{R}$ satisfying
$$\int_{\mathbb{R}^N}\Phi(|u|)dx<+\infty.$$
$L^{\Phi}(\mathbb{R}^N)$ is a Banach space endowed with Luxemburg norm
$$\|u\|_{\Phi}=\|u\|_{L^{\Phi}(\mathbb{R}^N)}:=\inf \left\{\lambda >0: \int_{\mathbb{R}^N}\Phi \left(\frac{|u|}{\lambda}\right)dx\leq1\right\}.$$
Particularly, when $\Phi(t)=|t|^p (p>1)$, the corresponding Orlicz space $L^{\Phi}(\mathbb{R}^N)$ reduces to the classical Lebesgue space $L^p(\mathbb{R}^N)$ endowed with the norm
$$\|u\|_p=L^{p}(\mathbb{R}^N):=\left(\int_{\mathbb{R}^N}|u(x)|^pdx\right)^{\frac{1}{p}}.$$
The fact that $\Phi$ satisfies $\Delta_2$-condition implies that
\begin{eqnarray}\label{ 2.2}
 u_n\rightarrow u \mbox{ in } L^{\Phi}(\Omega) \Longleftrightarrow \int_{\Omega}\Phi(|u_n-u|)dx\rightarrow 0.
 \end{eqnarray}
where $\Omega$ is an open set of $\mathbb{R}^N$. By the above Young's inequality (\ref{ 2.11}), a generalized type of H\"{o}lder's inequality (see \cite{Adams2003,M. M. Rao2002})
\begin{eqnarray}\label{ 2.3}
\left| \int_{\mathbb{R}^N}uvdx \right|\leq 2\|u\|_{\Phi}\|v\|_{\widetilde{\Phi}}, \quad \mbox{ for all } u \in L^{\Phi}(\mathbb{R}^N) \mbox{ and   } \mbox{all } v \in L^{\widetilde{\Phi}}(\mathbb{R}^N)
 \end{eqnarray}
can be obtained.
\par
Given an $N$-function $\Phi$ and a fractional parameter $0<s<1$, we recall the fractional Orlicz Sobolev space $W^{s, \Phi}(\mathbb{R}^N)$ defined as
$$W^{s, \Phi}(\mathbb{R}^N):=\left\{u \in L^{\Phi}(\mathbb{R}^N): \iint_{\mathbb{R}^{2N}}\Phi(|D_su|)d\mu<+\infty\right\},$$
where $D_su$ is defined by (\ref{1.3}) and $d\mu(x,y):=\frac{dxdy}{|x-y|^{N}}$. $W^{s, \Phi}(\mathbb{R}^N)$ is a Banach space endowed with the following norm
$$\|u\|_{s, \Phi}=\|u\|_{W^{s, \Phi}(\mathbb{R}^N)}:=\|u\|_{\Phi}+[u]_{s, \Phi},$$
where the so-called $(s, \Phi)$-Gagliardo semi-norm is defined as
$$[u]_{s, \Phi}:=\inf \left\{\lambda >0: \iint_{\mathbb{R}^{2N}}\Phi \left(\frac{|D_su|}{\lambda}\right)d\mu\leq1\right\}.$$
\par
The following lemmas will be useful in the sequel.
\vskip2mm
 \noindent
 {\bf Lemma 2.1.} (see \cite{Adams2003,Fukagai2006})\ \  {\it If $\Phi$ is an $N$-function, then the following conditions are equivalent:\\
$1)$
\begin{eqnarray}\label{ 2.4+}
1<l=\inf_{t>0}\frac{t\phi(t)}{\Phi(t)}\leq\sup_{t>0}\frac{t\phi(t)}{\Phi(t)}=m<+\infty;
\end{eqnarray}
$2)$ let $\zeta_0(t)=\min\{t^l, t^m\}$, $\zeta_1(t)=\max\{t^l, t^m\}$, for $t\geq0$. $\Phi$ satisfies
$$\zeta_0(t)\Phi(\rho)\leq \Phi(\rho t)\leq \zeta_1(t)\Phi(\rho), \quad \mbox{ for all } \rho, t\geq 0;$$
$3)$ $\Phi$ satisfies the $\Delta_2$-condition.
}
\vskip2mm
\noindent
 {\bf Lemma 2.2.} (see \cite{Bahrouni2020,Fukagai2006})\ \  {\it  If $\Phi$ is an $N$-function and (\ref{ 2.4+}) holds, then $\Phi$ satisfies \\}
$1)$
$$\zeta_0(\|u\|_{\Phi})\leq\int_{\mathbb{R}^N}\Phi(|u|)dx\leq\zeta_1(\|u\|_{\Phi}), \quad \mbox{ \it for all }u \in L^{\Phi}(\mathbb{R}^N);$$
$2)$
$$\zeta_0([u]_{s, \Phi})\leq\iint_{\mathbb{R}^{2N}}\Phi(|D_su|)d\mu\leq\zeta_1([u]_{s, \Phi}), \quad \mbox{ \it for all }u \in W^{s, \Phi}(\mathbb{R}^N).$$
 \vskip2mm
\noindent
 {\bf Lemma 2.3.} (see \cite{Fukagai2006})\ \  {\it If $\Phi$ is an $N$-function and (\ref{ 2.4+}) holds with $l>1$. Let $\widetilde{\Phi}$ be the complement of $\Phi$ and $\zeta_2(t)=\min\{t^{\widetilde{l}},t^{\widetilde{m}}\}$, $\zeta_3(t)=\max\{t^{\widetilde{l}},t^{\widetilde{m}}\}$, for $t\geq0$, where $\widetilde{l}:=\frac{l}{l-1}$ and $\widetilde{m}:=\frac{m}{m-1}$. Then $\widetilde{\Phi}$ satisfies\\}
$1)$
$$\widetilde{m}=\inf_{t>0}\frac{t\widetilde{\Phi}^{'}(t)}{\widetilde{\Phi}(t)}\leq\sup_{t>0}
\frac{t\widetilde{\Phi}^{'}(t)}{\widetilde{\Phi}(t)}=\widetilde{l};$$
$2)$
$$\zeta_2(t)\widetilde{\Phi}(\rho)\leq \widetilde{\Phi}(\rho t)\leq \zeta_3(t)\widetilde{\Phi}(\rho), \quad \mbox{\it for all } \rho, t\geq 0;$$
$3)$
$$\zeta_2(\|u\|_{\widetilde{\Phi}})\leq\int_{\mathbb{R}^N}\widetilde{\Phi}(|u|)dx\leq\zeta_3(\|u\|_{\widetilde{\Phi}}), \quad \mbox{\it for all }u \in L^{\widetilde{\Phi}}(\mathbb{R}^N).$$
\vskip2mm
\noindent
 {\bf Remark 2.4.}\ \   By Lemma 2.1 and Lemma 2.3, $(\phi_1)$ implies that $\Phi$ and $\widetilde{\Phi}$ are two $N$-functions satisfying $\Delta_2$-condition. The fact that $\Phi$ and $\widetilde{\Phi}$ satisfy $\Delta_2$-condition implies that $L^{\Phi}(\mathbb{R}^N)$ and $W^{s, \Phi}(\mathbb{R}^N)$ are separable and reflexive Banach spaces. Moreover, $C_c^{\infty}(\mathbb{R}^N)$ is dense in $W^{s, \Phi}(\mathbb{R}^N)$ (see \cite{Bonder2019,Adams2003,M. M. Rao2002}).
 \par
Next, we recall the Sobolev conjugate function of $\Phi$, which is denoted by $\Phi_*$. Suppose that
\begin{eqnarray}\label{2.6+}
 \int_0^1\frac{\Phi^{-1}(\tau)}{\tau^{\frac{N+s}{N}}}d\tau<+\infty\quad\mbox{ and }\quad \int_1^{+\infty}\frac{\Phi^{-1}(\tau)}{\tau^{\frac{N+s}{N}}}d\tau=+\infty.
 \end{eqnarray}
Then, $\Phi_*$ is defined by
$$\Phi_{*}^{-1}(t)=\int_{0}^{t}\frac{\Phi^{-1}(\tau)}{\tau^{\frac{N+s}{N}}}d\tau,\quad \mbox{ \it for } t\geq 0.$$
\vskip2mm
\noindent
 {\bf Lemma 2.5.} (see \cite{Bonder2023,Bahrouni2023})\ \  {\it If $\Phi$ is an $N$-function and (\ref{ 2.4+}) holds with $l,m\in(1, \frac{N}{s})$, then $(\ref{2.6+})$ holds. Let $\zeta_4(t)=\min\{t^{l^*},t^{m^*}\}$, $\zeta_5(t)=\max\{t^{l^*},t^{m^*}\}$, for $t\geq0$, where $l^*:=\frac{Nl}{N-sl}$, $m^*:=\frac{Nm}{N-sm}$. Then $\Phi_{*}$ satisfies}\\
$1)$
$$l^*=\inf_{t>0}\frac{t\Phi_{*}'(t)}{\Phi_*(t)}\leq\sup_{t>0}\frac{t\Phi_{*}'(t)}{\Phi_*(t)}=m^*;$$
$2)$
$$\zeta_4(t)\Phi_*(\rho)\leq \Phi_*(\rho t)\leq \zeta_5(t)\Phi_*(\rho), \quad \mbox{\it for all } \rho, t\geq 0;$$
$3)$
$$\zeta_4(\|u\|_{\Phi_*})\leq\int_{\mathbb{R}^N}\Phi_*(|u|)dx\leq\zeta_5(\|u\|_{\Phi_*}), \quad \mbox{\it for all }u \in L^{\Phi_*}(\mathbb{R}^N).$$
\par
The conjugate function $\Phi_*$ plays a crucial role in the following embedding results, which will be used frequently in our proofs.
\vskip2mm
\noindent
 {\bf Lemma 2.6.} (see \cite{Silva2021,Bahrouni2023})\ \  {\it If $\Phi$ is an $N$-function and (\ref{ 2.4+}) holds with $l,m\in(1, \frac{N}{s})$, then the following embedding results hold:\\
 $1)$ the embedding $W^{s, \Phi}(\mathbb{R}^N)\hookrightarrow L^{\Phi_{*}}(\mathbb{R}^N)$ is continuous;\\
 $2)$ the embedding $W^{s, \Phi}(\mathbb{R}^N)\hookrightarrow L^{\Phi}(\mathbb{R}^N)$ is continuous;\\
 $3)$ the embedding $W^{s, \Phi}(\mathbb{R}^N)\hookrightarrow L^\Psi(\mathbb{R}^N)$ is continuous for any $N$-function $\Psi$ satisfying $\Delta_2$-condition, $\Psi\prec\prec\Phi_{*}$ and
 $$\lim_{t\rightarrow 0^+}\frac{\Psi(t)}{\Phi(t)}=0;$$
 $4)$ when $\mathbb{R}^N$ is replaced by a $C^{0, 1}$ bounded open subset $D$ of $\mathbb{R}^N$, then the embedding $W^{s, \Phi}(D)\hookrightarrow L^\Psi(D)$ is compact for any $N$-function $\Psi$ satisfying $\Psi\prec\prec\Phi_{*}$.} Explicitly, when $m<l^*$, the embedding $W^{s, \Phi}(B_r)\hookrightarrow L^\Phi(B_r)$ is compact, where and in the sequel $B_r:=\{x\in \mathbb{R}^N: |x|<r\}$ for $r>0$.
 \vskip2mm
 \noindent
 {\bf Notation:}  Throughout this paper, $C_d$ is used to denote a positive constant which depends only on the constant or function $d$.
 \vskip-2mm
 \noindent
 \section{Proofs}
 \setcounter{equation}{0}
 \noindent
  \par
  On fractional Orlicz Sobolev space $W^{s, \Phi}(\mathbb{R}^N)$, denoted by $W$ for simplicity, the energy functional $I$ associated to equation (\ref{eq1}) is defined by
\begin{eqnarray}\label{3.1}
 I(u):=\iint_{\mathbb{R}^{2N}}\Phi(|D_su|)d\mu+\int_{\mathbb{R}^N}V(x)\Phi(|u|)dx-\int_{\mathbb{R}^N}F(x,u)dx.
 \end{eqnarray}
 It follows $(f_1)$ that for any given constant $\varepsilon>0$, there exists a constant $C_\varepsilon>0$ such that
 \begin{eqnarray}\label{3.2}
 |f(x,t)|\leq \varepsilon \phi(|t|)+C_\varepsilon\Phi_*'(|t|) \mbox{ and }
 |F(x,t)|\leq \varepsilon\Phi(|t|)+C_\varepsilon\Phi_*(|t|), \mbox{ for all } (x,t)\in\mathbb{R}^N\times\mathbb{R}.
 \end{eqnarray}
 Thus, by using standard arguments as \cite{Salort2020}, we have that $I\in C^1(W,\mathbb{R})$ and its derivative is given by
 \begin{eqnarray}\label{3.3}
  \langle I'(u),v\rangle=\iint_{\mathbb{R}^{2N}}a(|D_su|)D_su D_svd\mu
                                                   +\int_{\mathbb{R}^N}V(x)a(|u|)uvdx
                                                 -\int_{\mathbb{R}^N}f(x,u)vdx, \mbox{ for all }u, v\in W.
 \end{eqnarray}
 Thus, the critical points of $I$ are weak solutions of equation (\ref{eq1}).
 \par
 Define $I_i(i=1,2): W\rightarrow \mathbb{R}$ by
\begin{eqnarray}\label{3.4}
I_1(u)=\iint_{\mathbb{R}^{2N}}\Phi(|D_su|)d\mu+\int_{\mathbb{R}^N}V(x)\Phi(|u|)dx
\end{eqnarray}
 and
 \begin{eqnarray}\label{3.5}
 I_2(u)=\int_{\mathbb{R}^N}F(x,u)dx.
 \end{eqnarray}
 Then
 $$I(u)=I_1(u)-I_2(u), \mbox{ for all }u, v\in W$$
 and
 \begin{eqnarray}
 \label{3.5+}\langle I_1'(u),v\rangle&=&\iint_{\mathbb{R}^{2N}}a(|D_su|)D_su D_svd\mu
                                                   +\int_{\mathbb{R}^N}V(x)a(|u|)uvdx, \mbox{ for all }u, v\in W,\\
 \label{3.6+}\langle I_2'(u),v\rangle&=&\int_{\mathbb{R}^N}f(x,u)vdx, \mbox{ for all }u, v\in W.
 \end{eqnarray}
 \vskip2mm
 \noindent
 {\bf Lemma 3.1.}\ \  {\it Assume that $(\phi_1)$, $(V)$ and $(f_1)$ hold. Then there exist two constants $\rho, \eta>0$ such that $I(u)\geq \eta$ for all $u\in W$ with $\|u\|_{s, \Phi}=\rho$.}
 \vskip2mm
 \noindent
 {\bf Proof.}\  When $\|u\|_{s, \Phi}=\|u\|_{\Phi}+[u]_{s, \Phi}\leq 1$, by (\ref{3.1}), $(V)$, (\ref{3.2}) with taking $\varepsilon<\alpha_1$, Lemma 2.2, $3)$ in Lemma 2.5 and $1)$ in Lemma 2.6, we have
 \begin{eqnarray*}
 I(u)&\geq & \iint_{\mathbb{R}^{2N}}\Phi(|D_su|)d\mu+\alpha_1\int_{\mathbb{R}^N}\Phi(|u|)dx
               -\int_{\mathbb{R}^N}|F(x,u)|dx\nonumber\\
     &\geq & \iint_{\mathbb{R}^{2N}}\Phi(|D_su|)d\mu+(\alpha_1-\varepsilon)\int_{\mathbb{R}^N}\Phi(|u|)dx
               -C_\varepsilon\int_{\mathbb{R}^N}\Phi_*(|u|)dx\nonumber\\
     &\geq & [u]_{s, \Phi}^{m}+(\alpha_1-\varepsilon)\|u\|_{\Phi}^{m}
             -C_{\varepsilon}\max\{\|u\|_{\Phi_{*}}^{l^*},\|u\|_{\Phi_{*}}^{m^*}\}\nonumber\\
     &\geq & \min\{1,\alpha_1-\varepsilon\}C_{m}\| u\|_{s,\Phi}^{m}
             -C_{\varepsilon}C_{\Phi_{*}}^{l^*}\|u\|_{s,\Phi}^{l^*}-C_{\varepsilon}C_{\Phi_{*}}^{m^*}\|u\|_{s,\Phi}^{m^*}.
 \end{eqnarray*}
 Note that $m<l^*\leq m^*$. It is easy to see that the above inequality implies that there exist positive constants $\rho$ and $\eta$ small enough such that $I(u)\geq \eta$ for all $u\in W$ with $\|u\|_{s, \Phi}=\rho$. \quad $\Box$
\vskip2mm
 \noindent
 {\bf Lemma 3.2.}\ \  {\it Assume that $(\phi_1)$, $(V)$, $(f_1)$ and $(f_2)$ (or $(f_4)$) hold. Then there exists a $u_0\in W$ such that $I(tu_0)\rightarrow -\infty$ as $t\rightarrow+\infty$.}
 \vskip2mm
 \noindent
 {\bf Proof.}\  For any given constant $M>\alpha_2$, by $(f_1)$ and $(f_2)$ (or combine $(f_4)$ with $2)$ in Lemma 2.1), there exists a constant $C_M>0$ such that
 \begin{eqnarray}\label{3.6}
 F(x,t)\geq M\Phi(|t|)-C_M, \mbox{ for all } (x,t)\in\mathbb{R}^N\times\mathbb{R}.
 \end{eqnarray}
 Now, choose $u_0\in C_c^{\infty}(B_r)\setminus\{{\bf 0}\}$ with $0\leq u_0(x)\leq 1$. Then $u_0\in W,$ and by (\ref{3.1}), $(V)$, (\ref{3.6}), $2)$ in Lemma 2.1 and the fact $F(x,0)=0$ for all $x\in\mathbb{R}^N$, when $t>0$ we have
\begin{eqnarray*}
I(tu_0)& = & \iint_{\mathbb{R}^{2N}}\Phi(|D_s (tu_0)|)d\mu+\int_{\mathbb{R}^N}V(x)\Phi(|tu_0|)dx
            -\int_{\mathbb{R}^N}F(x,tu_0)dx\\
       & = & \iint_{\mathbb{R}^{2N}}\Phi(t|D_su_0|)d\mu+\int_{\mathbb{R}^N}V(x)\Phi(t|u_0|)dx
            -\int_{B_r}F(x,tu_0)dx\\
       &\leq&\Phi(t)\iint_{\mathbb{R}^{2N}}\max\{|D_su_0|^{l},|D_s u_0|^{m}\}d\mu
            +\alpha_2\int_{\mathbb{R}^N}\Phi(t|u_0|)dx-M\int_{B_r}\Phi(t|u_0|)+C_M|B_r|\\
       &\leq&\Phi(t)\iint_{\mathbb{R}^{2N}}(|D_su_0|^{l}+|D_s u_0|^{m})d\mu
            -(M-\alpha_2)\Phi(t)\int_{\mathbb{R}^N}\min\{|u_0|^{l},|u_0|^{m}\}dx+C_M|B_r|\\
       &\leq&\Phi(t)\left[\iint_{\mathbb{R}^{2N}}(|D_su_0|^{l}+|D_s u_0|^{m})d\mu
            -(M-\alpha_2)\|u_0\|_{m}^{m}\right]+C_M|B_r|.
 \end{eqnarray*}
 Note that $\lim\limits_{t\rightarrow+\infty}\Phi(t)=+\infty$. We can choose $M>\frac{\iint_{\mathbb{R}^{2N}}(|D_su_0|^{l}+|D_s u_0|^{m})d\mu}{\|u_0\|_{m_1}^{m_1}}+\alpha_2$ such that $I(tu_0)\rightarrow-\infty$ as $t\rightarrow+\infty$. What needs to be pointed out is that here we used the fact that $u_0\in W^{s, \Psi}(\mathbb{R}^N)$, where $\Psi(t)=|t|^l+|t|^m, t\geq0.$ So, $\iint_{\mathbb{R}^{2N}}(|D_su_0|^{l}+|D_s u_0|^{m})d\mu<+\infty$. \quad$\Box$
 \par
 Lemma 3.1, Lemma 3.2 and the fact $I({\bf 0})=0$ show that the energy functional $I$ has a mountain pass geometry: that is, setting
 $$\Gamma=\{\gamma\in C([0,1],W):\gamma(0)={\bf 0}, \|\gamma(1)\|_{s,\Phi}>\rho \mbox{ and } I(\gamma(1))\leq0\},$$
 we have $\Gamma\neq\emptyset$. Then, by using the variant version of the mountain pass theorem (see Theorem 3 in \cite{Silva2010}), we deduce that $I$ possesses a $(C)_c$-sequence $\{u_n\}$ with the level $c\geq\eta>0$ given by
 \begin{eqnarray}\label{3.7}
 c=\inf_{\gamma\in\Gamma}\max_{t\in [0,1]}I(\gamma(t)).
 \end{eqnarray}
 We recall that $(C)_c$-sequence $\{u_n\}$ of $I$ in $W$ means
  \begin{eqnarray}\label{3.8}
  I(u_n)\rightarrow c \quad  \mbox{ and } \quad (1+\|u_n\|_{s, \Phi})\|I'(u_n)\|_{W^*}\rightarrow 0, \quad \mbox{ as } n\rightarrow \infty.
 \end{eqnarray}
 \par
 To prove the boundedness of the $(C)_c$-sequence $\{u_n\}$ of $I$ in $W$, we will use the Lions' type result for fractional Orlicz-Sobolev spaces (see Theorem 1.6 in \cite{Silva2021}). We note that the claim $u_n\rightharpoonup 0$ in $X$ of Theorem 1.6 in \cite{Silva2021} is not necessary. With the same proof as Theorem 1.6 in \cite{Silva2021}, we can get the following result.
  \vskip2mm
 \noindent
 {\bf Lemma 3.3\it (Lions' type result for fractional Orlicz-Sobolev spaces).}\ \  {\it Suppose that the function $\phi$ defined by (\ref{1.2}) satisfies $(\phi_1)$ and
  $$\lim_{t\rightarrow 0^+}\frac{\Psi(t)}{\Phi(t)}=0.$$
  Let $\{u_n\}$ be a bounded sequence in $W^{s, \Phi}(\mathbb{R}^N)$ in such way that
  $$\lim_{n\rightarrow \infty}\sup_{y\in\mathbb{R}^N}\int_{B_r(y)}\Phi(|u_n|)dx=0,$$
  for some $r>0$. Then, $u_n\rightarrow \bf{0}$ in $L^{\Psi}(\mathbb{R}^N)$, where $\Psi$ is an $N$-function such that $\Psi\prec\prec\Phi_{*}$.
 }
  \vskip2mm
 \noindent
 {\bf Lemma 3.4.}\ \  {\it Assume that $(\phi_1)$, $(\phi_2)$, $(V)$ and $(f_1)$-$(f_3)$ hold. Then any $(C)_c$-sequence of $I$ in $W$ is bounded for all $c\geq0$.
 }
 \vskip2mm
 \noindent
 {\bf Proof.}\  Let $\{u_n\}$ be a $(C)_c$-sequence of $I$ in $W$ for $c\geq0$. By (\ref{3.8}), we have
 \begin{eqnarray}\label{3.9}
  I(u_n)\rightarrow c \mbox{ and }
  \left|\left\langle I'(u_n), \frac{1}{m}u_n\right\rangle\right|\rightarrow 0,  \mbox{ as } n\rightarrow \infty.
 \end{eqnarray}
 Then, by (\ref{3.1}), (\ref{3.3}), $(\phi_1)$ and $(V)$, for $n$ large, we have
  \begin{eqnarray}\label{3.10}
  c+1&\geq& I(u_n)-\left\langle I'(u_n), \frac{1}{m}u_n\right\rangle\nonumber\\
     &  = & \iint_{\mathbb{R}^{2N}}\left(\Phi(|D_su_n|)-\frac{1}{m}a(|D_su_n|)|D_su_n|^2\right)d\mu\nonumber\\
     &    & +\int_{\mathbb{R}^N}V(x)\left(\Phi(|u_n|)-\frac{1}{m}a(|u_n|)u_n^2\right)dx\nonumber\\
     &    & + \int_{\mathbb{R}^N}\left(\frac{1}{m}u_nf(x,u_n)-F(x,u_n)\right)dx\nonumber\\
     &\geq& \frac{1}{m}\int_{\mathbb{R}^N}\widehat{F}(x,u_n)dx.
 \end{eqnarray}
 \par
 To prove the boundedness of $\{u_n\}$, arguing by contradiction, we suppose that there exists a subsequence of $\{u_n\}$, still denoted by $\{u_n\}$, such that $\|u_n\|_{s, \Phi}\rightarrow \infty$, as $n\rightarrow \infty$. Let $\tilde{u}_n=\frac{u_n}{\|u_n\|_{s,\Phi}}$. Then $\|{\tilde{u}_n}\|_{s,\Phi}=1$.
 \par
  Firstly, we claim that
 \begin{eqnarray}\label{3.11}
 \lambda_1:=\lim_{n\rightarrow \infty}\sup_{y\in\mathbb{R}^N}\int_{B_2(y)}\Phi(|\tilde{u}_n|)dx=0.
 \end{eqnarray}
 Indeed, if $\lambda_1\neq 0$, there exist a constant $\delta>0$, a subsequence of $\{\tilde{u}_n\}$, still denoted by $\{\tilde{u}_n\}$, and a sequence $\{z_n\}\in \mathbb{Z}^N$ such that
 \begin{eqnarray}\label{3.12}
 \int_{B_2(z_n)}\Phi(|\tilde{u}_n|)dx>\delta, \quad \mbox{ for all } n\in \mathbb{N}.
 \end{eqnarray}
 Let $\bar{u}_n=\tilde{u}_n(\cdot +z_n)$. Then $\|\bar{u}_n\|_{s,\Phi}=\|\tilde{u}_n\|_{s,\Phi}=1$, that is, $\{{\bar{u}_n}\}$ is bounded in $W$. Passing to a subsequence of $\{{\bar{u}_n}\}$, still denoted by $\{{\bar{u}_n}\}$, by Remark 2.4 and $4)$ in Lemma 2.6, we can assume that there exists a $\bar{u}\in W$ such that
 \begin{eqnarray}\label{3.13}
 \bar{u}_n\rightharpoonup \bar{u} \mbox{ in } W, \quad \bar{u}_n\rightarrow \bar{u} \mbox{ in } L^{\Phi}(B_2) \quad \mbox{ and } \quad\bar{u}_n(x)\rightarrow \bar{u}(x) \mbox{ a.e. in } B_2.
 \end{eqnarray}
 Note that
 \begin{eqnarray*}
 \int_{B_2}\Phi(|\bar{u}_n|)dx=\int_{B_2(z_n)}\Phi(|\tilde{u}_n|)dx.
 \end{eqnarray*}
 Then, by (\ref{3.12}), (\ref{3.13}) and (\ref{ 2.2}), we obtain that $\bar{u}\neq{\bf 0}$ in $L^{\Phi}(B_2)$, that is, $[\bar{u}\neq0]:=\{x \in B_2: \bar u(x)\neq0\}$ has nonzero Lebesgue measure. Let $u_n^*=u_n(\cdot +z_n)$. Then $\|u_n^*\|_{s, \Phi}=\|u_n\|_{s, \Phi}$, and it follows from the fact that $V$ and $f$ are $1$-periodic in $x$ that
 \begin{eqnarray*}
I(u_n^*)=I(u_n) \quad\mbox{ and } \quad \|I'(u_n^*)\|_{W^*}=\|I'(u_n)\|_{W^*}, \quad \mbox{ for all } n\in\mathbb{N},
 \end{eqnarray*}
 which imply that $\{u_n^*\}$ is also a $(C)_c$-sequence of $I$. Then, by (\ref{3.10}), for $n$ large, we have
\begin{eqnarray}\label{3.14}
 \int_{\mathbb{R}^N}\widehat{F}(x,u_n^*)dx\leq m(c+1).
 \end{eqnarray}
 However, by $2)$ in Lemma 2.1, $(f_2)$ and $(f_3)$ imply
 \begin{eqnarray}\label{3.15}
 \lim_{|t|\rightarrow\infty}\widehat{F}(x,t)=+\infty, \quad \mbox{ uniformly in } x\in\mathbb{R}^N.
  \end{eqnarray}
 Moreover, by (\ref{3.13}), $\bar{u}_n=\tilde{u}_n(\cdot +z_n)=\frac{u_n(\cdot +z_n)}{\|u_n\|_{s,\Phi}}=\frac{u_n^*}{\|u_n\|_{s,\Phi}}$ implies
 \begin{eqnarray}\label{3.16}
 |u_n^*(x)|=|\bar{u}_n(x)|\|u_n\|_{s,\Phi}\rightarrow \infty, \quad \mbox{ a.e. } x\in [\bar{u}\neq0].
 \end{eqnarray}
  Then, it follows from $(f_3)$, (\ref{3.15}), (\ref{3.16}) and Fatou's Lemma that
 \begin{eqnarray*}
 \int_{\mathbb{R}^N}\widehat{F}(x,u_n^*)dx\geq \int_{[\bar{u}\neq0]}\widehat{F}(x,u_n^*)dx\rightarrow +\infty,\quad \mbox{ as } n\rightarrow\infty,
 \end{eqnarray*}
 which contradicts (\ref{3.14}). Therefore, $\lambda_1=0$ and thus (\ref{3.11}) holds.
 \par
 Next, for given $p\in[l,l^*)$ and $c>0$, by $(\phi_1)$, $(\phi_2)$ and $2)$ in Lemma 2.5, we have
 \begin{eqnarray}\label{3.17}
  \lim_{t\rightarrow 0^+}\frac{t^{p}}{\Phi(t)}=0\quad \mbox{ and} \quad \lim_{t\rightarrow+\infty}\frac{(ct)^{p}}{\Phi_{*}(t)}\leq\lim_{t\rightarrow
  +\infty}\frac{c^pt^{p}}{\Phi_{*}(1)\min\{t^{l^*},t^{m^*}\}}=0.
  \end{eqnarray}
 Then, by Lemma 3.3, (\ref{3.11}) and (\ref{3.17}) imply that
 \begin{eqnarray}\label{3.18}
 \tilde{u}_n\rightarrow {\bf 0} \mbox{ in } L^{p}(\mathbb{R}^N),\quad \mbox{ for all } p\in[l,l^*).
 \end{eqnarray}
 Hence, there exists a constant $M_1>0$ such that
 \begin{eqnarray}\label{3.19}
 \|\tilde{u}_n\|_{l}^{l}\leq M_1,\quad \mbox{ for all } n\in\mathbb{N}.
 \end{eqnarray}
 \par
 Finally, to get a contradiction, we will divide both sides of formula $I(u_n)=I_1(u_n)-I_2(u_n)$ by $\|u_n\|_{s,\Phi_1}^{l}$ and let $n\rightarrow\infty$. On one hand, by (\ref{3.9}), it is clear that
 \begin{eqnarray}\label{3.20}
\frac{I(u_n)}{\|u_n\|_{s,\Phi}^{l}}\rightarrow 0,\quad \mbox{ as } n\rightarrow\infty.
 \end{eqnarray}
 On the other hand, by (\ref{3.4}), $(V)$ and Lemma 2.2, we have
 \begin{eqnarray}\label{3.21}
 \frac{I_1(u_n)}{\|u_n\|_{s,\Phi}^{l}}
 & =  &\frac{\iint_{\mathbb{R}^{2N}}\Phi(|D_su_n|)d\mu+\int_{\mathbb{R}^N}V(x)\Phi(|u_n|)dx}
       {\|u_n\|_{s,\Phi}^{l}}\nonumber\\
 &\geq& \frac{\min\{[u_n]_{s, \Phi}^{l},[u_n]_{s,\Phi}^{m}\}
        +\alpha_1\min\{\|u_n\|_{\Phi}^{l},\|u_n\|_{\Phi}^{m}\}}{\|u_n\|_{s,\Phi}^{l}}\nonumber\\
 &\geq& \frac{[u_n]_{s, \Phi}^{l}+\alpha_1\|u_n\|_{\Phi}^{l}-1-\alpha_1}{\|u_n\|_{s,\Phi}^{l}}\nonumber\\
 &\geq& \frac{\min\{1,\alpha_1\}C_{l}([u_n]_{s, \Phi}+\|u_n\|_{\Phi})^{l}-1-\alpha_1}
        {\|u_n\|_{s,\Phi}^{l}}\rightarrow \min\{1,\alpha_1\}C_{l}, \mbox{ as } n\rightarrow\infty.
 \end{eqnarray}
 Moreover, by $2)$ in Lemma 2.1, $(f_1)$ implies that
 \begin{eqnarray*}
 \lim_{|t|\rightarrow 0}\frac{F(x,t)}{|t|^{l}}=0, \quad \mbox{ uniformly in } x\in \mathbb{R}^N.
 \end{eqnarray*}
 Then, for any given constant $\varepsilon>0$, there exists a constant $R_\varepsilon>0$ such that
\begin{eqnarray}\label{3.22}
 \frac{|F(x,t)|}{|u|^{l}}\leq\varepsilon,\quad \mbox{ for all } x\in\mathbb{R}^N, |t|\leq R_\varepsilon.
 \end{eqnarray}
 For above $R_\varepsilon>0$, by $(f_1)$ and $(f_3)$, there exists a constant $C_R>0$ such that
 \begin{eqnarray}\label{3.23}
 \left(\frac{|F(x,t)|}{|t|^{l}}\right)^k\leq C_R\widehat{F}(x,t),\quad \mbox{ for all } x\in\mathbb{R}^N,  |t|> R_\varepsilon.
 \end{eqnarray}
 Let
 $$X_n=\{x\in\mathbb{R}^N: |u_n(x)|\leq R_\varepsilon \} \mbox{ and } Y_n=\{x\in\mathbb{R}^N: |u_n(x)|> R_\varepsilon\}.$$
 Then
 \begin{eqnarray}\label{3.24}
 \frac{|I_2(u_n)|}{\|u_n\|_{s,\Phi}^{l}}
 \leq\int_{X_n}\frac{|F(x,u_n)|}{\|u_n\|_{s,\Phi}^{l}}dx+\int_{Y_n}\frac{|F(x,u_n)|}{\|u_n\|_{s,\Phi}^{l}}dx.
 \end{eqnarray}
 By (\ref{3.22}) and (\ref{3.19}), we have
 \begin{eqnarray}\label{3.25}
 \int_{X_n}\frac{|F(x,u_n)|}{\|u_n\|_{s,\Phi}^{l}}dx
 =\int_{X_n}\frac{|F(x,u_n)|}{|u_n|^{l}}|\tilde{u}_n|^{l}dx\leq \varepsilon\|\tilde{u}_n\|_{l}^{l}
 \leq \varepsilon M_1.
 \end{eqnarray}
  The claim $k>\frac{N}{sl}$ given by $(f_3)$ implies that $\frac{lk}{k-1}\in(l,l^*)$. Hence, by H\"{o}lder's inequality, (\ref{3.23}), (\ref{3.10}), (\ref{3.18}) and the fact $\widehat{F}(x,t)\geq 0$, we have
 \begin{eqnarray}\label{3.26}
 \int_{Y_n}\frac{|F(x,u_n)|}{\|u_n\|_{s,\Phi}^{l}}dx
 &=&\int_{Y_n}\frac{|F(x,u_n)|}{|u_n|^{l}}|\tilde{u}_n|^{l}dx\nonumber\\
 &\leq& \left(\int_{Y_n}\left(\frac{|F(x,u_n)|}{|u_n|^{l}}\right)^kdx\right)^{\frac{1}{k}}
        \left(\int_{Y_n}|\tilde{u}_n|^{\frac{lk}{k-1}}dx\right)^{\frac{k-1}{k}}\nonumber\\
  &\leq& \left(\int_{Y_n}C_R\widehat{F}(x,u_n)dx\right)^{\frac{1}{k}}
         \|\tilde{u}_n\|_{\frac{lk}{k-1}}^l\nonumber\\
  &\leq& [C_Rm(c+1)]^{\frac{1}{k}}\|\tilde{u}_n\|_{\frac{lk}{k-1}}^l \rightarrow 0,\quad \mbox{ as } n\rightarrow\infty.
 \end{eqnarray}
 Since $\varepsilon$ is arbitrary, it follows from (\ref{3.25}) and (\ref{3.26}) that
 \begin{eqnarray}\label{3.27}
 \frac{I_2(u_n)}{\|u_n\|_{s,\Phi}^{l}}\rightarrow 0,\quad \mbox{ as } n\rightarrow\infty.
 \end{eqnarray}
 By dividing both sides of formula $I(u_n)=I_1(u_n)-I_2(u_n)$ by $\|u_n\|_{s,\Phi_1}^{l}$ and letting $n\rightarrow\infty$, we get a contradiction via (\ref{3.20}), (\ref{3.21}) and (\ref{3.27}). Therefore, the $(C)_c$-sequence $\{u_n\}$ is bounded. \quad $\Box$
  \vskip2mm
 \noindent
 {\bf Lemma 3.5.}\ \  {\it Assume that $(\phi_1)$, $(V)$, $(f_1)$, $(f_4)$ and $(f_5)$ are satisfied. Then for $u\in W$, there holds
 $$I(u)\geq I(tu)+\frac{1-t^l}{m}\langle I'(u),u\rangle, \quad \mbox{ for all } t\in[0, \theta_0],$$
 where $\theta_0$ is given in $(f_5)$.
 }
 \vskip2mm
 \noindent
 {\bf Proof.}\  When $u\in W$, $0\leq t\leq 1$, by (\ref{3.1}), (\ref{3.3}) and Lemma 2.1, we have
  \begin{eqnarray*}
 & &I(u)- I(tu)-\frac{1-t^l}{m}\langle I'(u),u\rangle\nonumber\\
 &=&\iint_{\mathbb{R}^{2N}}\Phi(|D_su|)d\mu+\int_{\mathbb{R}^N}V(x)\Phi(|u|)dx
    -\int_{\mathbb{R}^N}F(x,u)dx\nonumber\\
 & &-\iint_{\mathbb{R}^{2N}}\Phi(|D_stu|)d\mu-\int_{\mathbb{R}^N}V(x)\Phi(|tu|)dx
    +\int_{\mathbb{R}^N}F(x,tu)dx\nonumber\\
 & &-\frac{1-t^l}{m}\iint_{\mathbb{R}^{2N}}a(|D_su|)|D_su|^2d\mu
    -\frac{1-t^l}{m}\int_{\mathbb{R}^N}V(x)a(|u|)u^2dx
    +\frac{1-t^l}{m}\int_{\mathbb{R}^N}uf(x,u)dx\nonumber\\
 &\geq&\iint_{\mathbb{R}^{2N}}\Phi(|D_su|)d\mu
       -\max\{t^l, t^m\}\iint_{\mathbb{R}^{2N}}\Phi(|D_su|)d\mu
       -(1-t^l)\iint_{\mathbb{R}^{2N}}\Phi(|D_su|)d\mu\nonumber\\
  &   &+\int_{\mathbb{R}^N}V(x)\Phi(|u|)dx
       -\max\{t^l, t^m\}\int_{\mathbb{R}^N}V(x)\Phi(|u|)dx
       -(1-t^l)\int_{\mathbb{R}^N}V(x)\Phi(|u|)dx\nonumber\\
  &   &+\int_{\mathbb{R}^N}\left[\frac{1-t^l}{m}uf(x,u)-F(x, u)+F(x, tu)\right]dx\nonumber\\
  & = &\int_{\mathbb{R}^N}\left[\frac{1-t^l}{m}uf(x,u)-\int_{tu}^{u}f(x,\tau)d\tau\right]dx.
 \end{eqnarray*}
  Then, it follows from $(f_5)$ that
  $$I(u)\geq I(tu)+\frac{1-t^l}{m}\langle I'(u),u\rangle, \quad \mbox{ for all } t\in[0, \theta_0],$$
 for some $\theta_0\in (0, 1)$. \quad $\Box$
 \vskip2mm
 \noindent
 {\bf Lemma 3.6.}\ \  {\it Assume that $(\phi_1)$, $(V)$, $(f_1)$, $(f_4)$ and $(f_5)$ hold. Then any $(C)_c$-sequence of $I$ in $W$ is bounded for all $c\geq0$.
 }
 \vskip2mm
 \noindent
 {\bf Proof.}\  Let $\{u_n\}$ be a $(C)_c$-sequence of $I$ in $W$ for $c\geq0$. By (\ref{3.8}), we have
 \begin{eqnarray}\label{3.28}
  I(u_n)\rightarrow c \mbox{ and }
  \left|\left\langle I'(u_n), u_n\right\rangle\right|\rightarrow 0,  \mbox{ as } n\rightarrow \infty.
 \end{eqnarray}
 To prove the boundedness of $\{u_n\}$, arguing by contradiction, we suppose that there exists a subsequence of $\{u_n\}$, still denoted by $\{u_n\}$, such that $\|u_n\|_{s, \Phi}\rightarrow \infty$, as $n\rightarrow \infty$. Let $\tilde{u}_n=\frac{u_n}{\|u_n\|_{s,\Phi}}$. Then $\|{\tilde{u}_n}\|_{s,\Phi}=1$.
 \par
  Firstly, we claim that
 \begin{eqnarray}\label{3.29}
 \lambda_2:=\lim_{n\rightarrow \infty}\sup_{y\in\mathbb{R}^N}\int_{B_2(y)}\Phi(|\tilde{u}_n|)dx=0.
 \end{eqnarray}
 Indeed, if $\lambda_2\neq 0$, there exist a constant $\delta>0$, a subsequence of $\{\tilde{u}_n\}$, still denoted by $\{\tilde{u}_n\}$, and a sequence $\{z_n\}\in \mathbb{Z}^N$ such that
 \begin{eqnarray}\label{3.30}
 \int_{B_2(z_n)}\Phi(|\tilde{u}_n|)dx>\delta, \quad \mbox{ for all } n\in \mathbb{N}.
 \end{eqnarray}
 Let $\bar{u}_n=\tilde{u}_n(\cdot +z_n)$. Then $\|\bar{u}_n\|_{s,\Phi}=\|\tilde{u}_n\|_{s,\Phi}=1$, that is, $\{{\bar{u}_n}\}$ is bounded in $W$. Passing to a subsequence of $\{{\bar{u}_n}\}$, still denoted by $\{{\bar{u}_n}\}$, by Remark 2.4 and $4)$ in Lemma 2.6, we can assume that there exists a $\bar{u}\in W$ such that
 \begin{eqnarray}\label{3.31}
 \bar{u}_n\rightharpoonup \bar{u} \mbox{ in } W, \quad \bar{u}_n\rightarrow \bar{u} \mbox{ in } L^{\Phi}(B_2) \quad \mbox{ and } \quad\bar{u}_n(x)\rightarrow \bar{u}(x) \mbox{ a.e. in } B_2.
 \end{eqnarray}
 Note that
 \begin{eqnarray*}
 \int_{B_2}\Phi(|\bar{u}_n|)dx=\int_{B_2(z_n)}\Phi(|\tilde{u}_n|)dx.
 \end{eqnarray*}
 Then, by (\ref{3.30}), (\ref{3.31}) and (\ref{ 2.2}), we obtain that $\bar{u}\neq{\bf 0}$ in $L^{\Phi}(B_2)$, that is, $[\bar{u}\neq0]:=\{x \in B_2: \bar u(x)\neq0\}$ has nonzero Lebesgue measure. Let $u_n^*=u_n(\cdot +z_n)$. Then $\|u_n^*\|_{s, \Phi}=\|u_n\|_{s, \Phi}$, and
 \begin{eqnarray}\label{3.32}
 |u_n^*(x)|=|\bar{u}_n(x)|\|u_n\|_{s,\Phi}\rightarrow \infty, \quad \mbox{ a.e. } x\in [\bar{u}\neq0].
 \end{eqnarray}
  Then, it follows from (\ref{3.5}), $(f_4)$, (\ref{3.32}) and Fatou's Lemma that
  \begin{eqnarray}\label{3.33}
 \frac{I_2(u_n)}{\|u_n\|_{s,\Phi}^{m}}
 &=&\int_{\mathbb{R}^N}\frac{F(x,u_n)}{\|u_n\|_{s,\Phi}^{m}}dx\nonumber\\
 &=&\int_{\mathbb{R}^N}\frac{F(x+z_n,u_n^*)}{|u_n^*|^{m}}|\bar{u}_n|^mdx\nonumber\\
 &\geq&\int_{[\bar{u}\neq0]}\frac{F(x+z_n,u_n^*)}{|u_n^*|^{m}}|\bar{u}_n|^mdx\rightarrow +\infty,\quad \mbox{ as } n\rightarrow\infty.
 \end{eqnarray}
 Moreover, it follows from (\ref{3.4}), $(V)$ and Lemma 2.2 that
 \begin{eqnarray}\label{3.34}
 \limsup_{n\rightarrow \infty}\frac{I_1(u_n)}{\|u_n\|_{s,\Phi}^{m}}
 &=&\limsup_{n\rightarrow \infty}\frac{\iint_{\mathbb{R}^2N}\Phi(|D_su|)d\mu
    +\int_{\mathbb{R}^N}V(x)\Phi(|u|)dx}{\|u_n\|_{s,\Phi}^{m}}\nonumber\\
 &\leq&\limsup_{n\rightarrow \infty}\frac{\max\{[u_n]_{s, \Phi}^{l},[u_n]_{s,\Phi}^{m}\}
        +\alpha_2\max\{\|u_n\|_{\Phi}^{l},\|u_n\|_{\Phi}^{m}\}}{\|u_n\|_{s,\Phi}^{m}}\nonumber\\
 &\leq& 1+\alpha_2.
 \end{eqnarray}
 By dividing both sides of formula $I(u_n)=I_1(u_n)-I_2(u_n)$ by $\|u_n\|_{s,\Phi_1}^{m}$ and letting $n\rightarrow\infty$, we get a contradiction via (\ref{3.28}), (\ref{3.33}) and (\ref{3.34}). Therefore, $\lambda_2=0$ and thus (\ref{3.29}) holds. Then, by using the Lions' type result for fractional Orlicz-Sobolev spaces, with the similar discussion as in Lemma 3.4, we have
 \begin{eqnarray}\label{3.35}
 \tilde{u}_n\rightarrow {\bf 0} \mbox{ in } L^{p}(\mathbb{R}^N),\quad \mbox{ for all } p\in(m,l^*).
 \end{eqnarray}
 \par
 Besides, it follows from $1)$ in Lemma 2.2, $3)$ in Lemma 2.5, $1)$-$2)$ in Lemma 2.6 and the fact $\|{\tilde{u}_n}\|_{s,\Phi}=1$ that there exists a constant $M_2>0$ such that
  \begin{eqnarray}\label{3.36}
 &&\int_{\mathbb{R}^N}\left(\Phi(|\tilde{u}_n|)+\Phi_*(|\tilde{u}_n|)\right)dx\nonumber\\
 &\leq& \max\left\{\|\tilde{u}_n\|_{\Phi}^{l},\|\tilde{u}_n\|_{\Phi}^{m}\right\}
 +\max\left\{\|\tilde{u}_n\|_{\Phi_*}^{l^*},\|\tilde{u}_n\|_{\Phi_*}^{m^*}\right\}\nonumber\\
 &\leq&M_2, \quad \mbox{ for all } n\in\mathbb{N}.
 \end{eqnarray}
  \par
 Next, for any given $R>1$, let $t_n=\frac{R}{\|u_n\|_{s, \Phi}}$. Since $\|u_n\|_{s, \Phi}\rightarrow \infty$ as $n\rightarrow \infty$, it follows that $t_n\in(0, \theta_0)$ for $n$ large enough. Thus, by (\ref{3.28}) and Lemma 3.5, we have
 \begin{eqnarray}\label{3.37}
 c+o_n(1)
 &=&I(u_n)\nonumber\\
 &\geq&I(t_nu_n)+\frac{1-t_n^l}{m}\langle I'(u_n),u_n\rangle\nonumber\\
 &=& I\left(\frac{R}{\|u_n\|_{s, \Phi}}u_n\right)+o_n(1)\nonumber\\
 &=& I(R\tilde{u}_n)+o_n(1)\nonumber\\
 &=& I_1(R\tilde{u}_n)-I_2(R\tilde{u}_n)+o_n(1).
 \end{eqnarray}
 For above $R$ and any given $\varepsilon>0$, by $(f_1)$, the continuity of $F$ and the fact that $\Phi$ and $\Phi_*$ satisfy the $\Delta_2$-condition, there exist constants $C_\varepsilon>0$ and $p\in(m, l^*)$ such that
 \begin{eqnarray}\label{3.38}
 |F(x,Rt)|\leq\varepsilon(\Phi(|t|)+\Phi_{*}(|t|))+C_\varepsilon|t|^{p}, \mbox{ for all } (x,t)\in\mathbb{R}^N\times\mathbb{R}.
 \end{eqnarray}
 Then, by (\ref{3.5}), (\ref{3.35}), (\ref{3.36}) and (\ref{3.38}), we have
 \begin{eqnarray}\label{3.39}
 |I_2(R\tilde{u}_n)|
 &\leq&\int_{\mathbb{R}^N}|F(x,R\tilde{u}_n)|dx\nonumber\\
 &\leq&\varepsilon\int_{\mathbb{R}^N}(\Phi(|\tilde{u}_n|)+\Phi_{*}(|\tilde{u}_n|))dx
 +C_\varepsilon\int_{\mathbb{R}^N}|\tilde{u}_n|^{p}dx\nonumber\\
 &\leq& \varepsilon M_2+o_n(1).
 \end{eqnarray}
 Since $\varepsilon>0$ is arbitrary, we conclude that
  \begin{eqnarray}\label{3.40}
 I_2(R\tilde{u}_n)=o_n(1).
 \end{eqnarray}
 Moreover, for above $R>1$, by (\ref{3.4}), Lemma 2.1 and the fact $\|{\tilde{u}_n}\|_{s,\Phi}=\|\tilde{u}_n\|_{\Phi}+[\tilde{u}_n]_{s, \Phi}=1$, we have
 \begin{eqnarray}\label{3.41}
I_1(R\tilde{u}_n)
& = &\iint_{\mathbb{R}^{2N}}\Phi(|D_s(R\tilde{u}_n)|)d\mu
    +\int_{\mathbb{R}^N}V(x)\Phi(|R\tilde{u}_n|)dx\nonumber\\
&\geq&\min\{R^{l},R^{m}\}\left(\min\{[\tilde{u}_n]_{s, \Phi}^{l},[\tilde{u}_n]_{s,\Phi}^{m}\}
      +\alpha_1\min\{\|\tilde{u}_n\|_{\Phi}^{l},\|\tilde{u}_n\|_{\Phi}^{m}\}\right)\nonumber\\
& = &R^{l}\left([\tilde{u}_n]_{s,\Phi}^{m}+\alpha_1\|\tilde{u}_n\|_{\Phi}^{m}\right)\nonumber\\
&\geq&\min\{1,\alpha_1\}R^{l}\left([\tilde{u}_n]_{s,\Phi}^{m}+\|\tilde{u}_n\|_{\Phi}^{m}\right)\nonumber\\
&\geq&\min\{1,\alpha_1\}R^{l}C_m,
\end{eqnarray}
where $C_m:=\inf_{|u|+|v|=1}\{|u|^m+|v|^m\}>0$. Then, by the arbitrariness of $R$, combining (\ref{3.40}) and (\ref{3.41}) with (\ref{3.37}), we get a contradiction. Therefore, the $(C)_c$-sequence $\{u_n\}$ is bounded. \quad $\Box$
\vskip2mm
 \noindent
 {\bf Lemma 3.7.}\ \  {\it Assume that $(\phi_1)$, $(V)$ and $(f_1)$ hold. Then $I': W\rightarrow W^*$ is weakly sequentially continuous. Namely, if $u_n\rightharpoonup u$ in $W$, then $I'(u_n)\rightharpoonup I'(u)$ in the dual space $W^*$ of $W$.
 }
 \vskip2mm
 \noindent
 {\bf Proof.}\  Since $W$ is reflexive, it is enough to prove $I'(u_n)\stackrel{w^*}{\rightharpoonup}I'(u)$ in $W^*$. Namely, to prove
 \begin{eqnarray}\label{3.42}
\lim_{n\rightarrow \infty}\langle I'(u_n),v\rangle=\langle I'(u),v\rangle, \mbox{ for all } v\in W.
  \end{eqnarray}
\par
Firstly, we prove that $I'$ is bounded on each bounded subset of $W$. Indeed, by (\ref{3.3}), $(V)$, (\ref{ 2.11}), (\ref{3.2}), (\ref{ 2.1}), Lemma 2.2, $3)$ in Lemma 2.5, $1)$ in Lemma 2.6 and the fact that $\Phi$, $\widetilde{\Phi}$ and $\Phi_*$ satisfy the $\Delta_2$-condition, we have
 \begin{eqnarray*}
\|I'(u)\|_{W^*}
&=&\sup_{v\in W, \|v\|_{s, \Phi}=1}|\langle I'(u),v\rangle|\nonumber\\
&\leq&\sup_{v\in W, \|v\|_{s, \Phi}=1}\left(\iint_{\mathbb{R}^{2N}}a(|D_su|)|D_su||D_sv|d\mu
      +\int_{\mathbb{R}^N}V(x)a(|u|)|u||v|dx\right.\nonumber\\
&&\left.+\int_{\mathbb{R}^N}|f(x,u)||v|dx\right)\nonumber\\
&\leq&\sup_{v\in W, \|v\|_{s, \Phi}=1}\left(\iint_{\mathbb{R}^{2N}}\widetilde{\Phi}(a(|D_su|)|D_su|)d\mu
      +\iint_{\mathbb{R}^{2N}}\Phi(|D_sv|)d\mu\right.\nonumber\\
&&+(\alpha_2+\varepsilon)\int_{\mathbb{R}^N}\widetilde{\Phi}(a(|u|)|u|)dx
  +(\alpha_2+\varepsilon)\int_{\mathbb{R}^N}\Phi(|v|)dx\nonumber\\
&&\left.+C_{\varepsilon}\int_{\mathbb{R}^N}\widetilde{\Phi}_*(\Phi_*'(|u|))dx
  +C_{\varepsilon}\int_{\mathbb{R}^N}\Phi_*(|v|)dx\right)\nonumber\\
&\leq&\left(\iint_{\mathbb{R}^{2N}}\Phi(2|D_su|)d\mu
      +(\alpha_2+\varepsilon)\int_{\mathbb{R}^N}\Phi(2|u|)dx
+C_{\varepsilon}\int_{\mathbb{R}^N}\Phi_*(2|u|)dx\right)\nonumber\\
&& +\sup_{v\in W, \|v\|_{s, \Phi}=1}\left(\max\{[v]_{s, \Phi}^{l},[v]_{s,\Phi}^{m}\}
        +(\alpha_2+\varepsilon)\max\{\|v\|_{\Phi}^{l},\|v\|_{\Phi}^{m}\}\right.\nonumber\\
&&       \left.+C_{\varepsilon}\max\{\|v\|_{\Phi_*}^{l^*},\|v\|_{\Phi_*}^{m^*}\}\right)\nonumber\\
&\leq&K_2\left(\iint_{\mathbb{R}^{2N}}\Phi(|D_su|)d\mu
      +(\alpha_2+\varepsilon)\int_{\mathbb{R}^N}\Phi(|u|)dx
+C_{\varepsilon}\int_{\mathbb{R}^N}\Phi_*(|u|)dx\right)\nonumber\\
&& +1+\alpha_2+\varepsilon+C_{\varepsilon}C_{\Phi_*}\nonumber\\
&\leq&K_2\left((1+\alpha_2+\varepsilon)\|u\|_{s,\Phi}^{m}
      +C_{\varepsilon}C_{\Phi_*}\|u\|_{s,\Phi}^{m^*}\right)+(K_2+1)(1+\alpha_2+\varepsilon+C_{\varepsilon}C_{\Phi_*}),
  \end{eqnarray*}
which implies that $I'$ is bounded on each bounded subset of $W$. Moreover, $C_c^{\infty}(\mathbb{R}^N)$ is dense in $W$. Then, to prove (\ref{3.42}) we only need to prove
\begin{eqnarray}\label{3.43}
\lim_{n\rightarrow \infty}\langle I'(u_n),w\rangle=\langle I'(u),w\rangle, \mbox{ for all } w\in C_c^{\infty}(\mathbb{R}^N).
  \end{eqnarray}
\par
To get (\ref{3.43}), arguing by contradiction, we suppose that there exist constant $\delta>0$, $w_0\in C_c^{\infty}(\mathbb{R}^N)$ with $\mbox{supp}\{w_0\}\subset B_r$ for some $r>0$, and a subsequence of $\{u_n\}$, still denoted by $\{u_n\}$, such that
\begin{eqnarray}\label{3.44}
|\langle I'(u_n),w_0\rangle-\langle I'(u),w_0\rangle|\geq\delta, \mbox{ for all } n\in \mathbb{R}^N.
  \end{eqnarray}
Since $u_n\rightharpoonup u$ in $W$, by $4)$ in Lemma 2.6, there exists a subsequence of $\{u_n\}$, still denoted by $\{u_n\}$, such that
 \begin{eqnarray*}
 {u}_n\rightarrow u \mbox{ in } L_{loc}^{\Phi}(\mathbb{R}^N), \quad {u}_n(x)\rightarrow u(x) \mbox{ a.e. in } \mathbb{R}^{N} \quad \mbox{ and } \quad D_s{u}_n\rightarrow D_su \mbox{ a.e. in } \mathbb{R}^{2N}.
 \end{eqnarray*}
\par
 Next, we claim that
 \begin{eqnarray}\label{3.45}
 \lim_{n\rightarrow\infty}\int_{\mathbb{R}^N} f(x,u_n)w_0dx=\int_{\mathbb{R}^N}f(x,u)w_0dx.
 \end{eqnarray}
 Indeed, it follows $(f_1)$ that for any given constant $\varepsilon>0$, there exists a constant $C_\varepsilon>0$ such that
 \begin{eqnarray*}
 |f(x,t)|\leq C_\varepsilon+\varepsilon\Phi_*'(|t|), \mbox{ for all } (x,t)\in\mathbb{R}^N\times\mathbb{R}.
 \end{eqnarray*}
 Then, by using standard arguments, we can obtain that the sequence $\{f(x,u_n)\}$ is bounded in $L^{\widetilde{\Phi}_*}(B_r)$. Moreover, $f(x,u_n)\rightarrow f(x,u)$ a.e. in $B_r$. Then, by applying Lemma 2.1 in \cite{Alves2014}, we get $f(x,u_n)\rightharpoonup f(x,u)$ in $L^{\widetilde{\Phi}_*}(B_r)$, and thus (\ref{3.45}) holds because $w_0\in L^{\Phi_*}(B_r)$.
 \par
 Similarly, we can get
 \begin{eqnarray}\label{3.46}
 \lim_{n\rightarrow\infty}\iint_{\mathbb{R}^{2N}}a(|D_su_n|)D_su_nD_sw_0d\mu
 =\iint_{\mathbb{R}^{2N}}a(|D_su|)D_suD_sw_0d\mu
 \end{eqnarray}
 and
  \begin{eqnarray}\label{3.47}
 \lim_{n\rightarrow\infty}\int_{\mathbb{R}^{N}}V(x)a(|u_n|)u_nw_0dx
 =\int_{\mathbb{R}^{N}}V(x)a(|u|)uw_0dx,
 \end{eqnarray}
 which based on the fact that the sequence $\{a(|D_su_n|)D_su_n\}$ is bounded in $L^{\widetilde{\Phi}}(\mathbb{R}^{2N}, d\mu)$, $a(|D_su_n|)D_su_n\rightarrow a(|D_su|)D_su$ a.e. in $\mathbb{R}^{2N}$, $D_sw_0\in L^{\Phi}(\mathbb{R}^{2N}, d\mu)$, and the sequence $\{V(x)a(|u_n|)u_n\}$ is bounded in $L^{\widetilde{\Phi}}(\mathbb{R}^{N})$, $V(x)a(|u_n|)u_n\rightarrow V(x)a(|u|)u$ a.e. in $\mathbb{R}^{N}$, $w_0\in L^{\Phi}(\mathbb{R}^{N})$, respectively.
 \par
 Therefore, combining (\ref{3.45})-(\ref{3.47}) with (\ref{3.3}), we can conclude that
 \begin{eqnarray*}
 \lim_{n\rightarrow\infty}|\langle I'(u_n),w_0\rangle-\langle I'(u),w_0\rangle|=0,
  \end{eqnarray*}
  which contradicts (\ref{3.44}). So, (\ref{3.43}) holds and the proof is completed. \quad $\Box$
 \vskip2mm
 \noindent
 {\bf Lemma 3.8.}\ \  {\it Equation (\ref{eq1}) has at least a nontrivial solution under the assumptions of Theorem 1.1 and Theorem 1.2, respectively.
 }
 \vskip2mm
 \noindent
 {\bf Proof.}\  Let $\{u_n\}$ be the $(C)_c$-sequence of $I$ in $W$ for the level $c>0$ given in (\ref{3.7}). Lemma 3.4 and Lemma 3.6 show that the sequence $\{u_n\}$ is bounded in $W$ under the assumptions of Theorem 1.1 and Theorem 1.2, respectively.
 \par
 First, we claim that
 \begin{eqnarray}\label{3.48}
 \lambda_3:=\lim_{n\rightarrow \infty}\sup_{y\in\mathbb{R}^N}\int_{B_2(y)}\Phi(|u_n|)dx>0.
 \end{eqnarray}
 Indeed, if $\lambda_3=0$, by using the Lions' type result for fractional Orlicz-Sobolev spaces again, we have
 \begin{eqnarray}\label{3.49}
 u_n\rightarrow {\bf 0} \mbox{ in } L^{p}(\mathbb{R}^N),\quad \mbox{ for all } p\in(m,l^*).
 \end{eqnarray}
  Given $p\in(m,l^*)$, by $(f_1)$, $(\phi_1)$ and the definition $F(x,t)=\int_{0}^{t}f(x,\tau)d\tau$, for any given constant $\varepsilon>0$, there exists a constant $C_\varepsilon>0$ such that
   \begin{eqnarray}\label{3.50}
 |F(x,t)|\leq \varepsilon(\Phi(|t|)+\Phi_*(|t|))+C_\varepsilon|t|^p, \mbox{ for all } (x,t)\in\mathbb{R}^N\times\mathbb{R}
 \end{eqnarray}
 and
 \begin{eqnarray}\label{3.51}
 |tf(x,t)|\leq \varepsilon(\Phi(|t|)+\Phi_*(|t|))+C_\varepsilon|t|^p, \mbox{ for all } (x,t)\in\mathbb{R}^N\times\mathbb{R}.
 \end{eqnarray}
 Then, it follows from (\ref{3.49})-(\ref{3.51}), $1)$ in Lemma 2.2, $3)$ in Lemma 2.5 and $1)$ in Lemma 2.6, the boundedness of $\{u_n\}$ and the arbitrariness of $\varepsilon$ that
  \begin{eqnarray}\label{3.52}
 \lim_{n\rightarrow\infty}\int_{\mathbb{R}^N}F(x,u_n)dx
 =\lim_{n\rightarrow\infty}\int_{\mathbb{R}^N}u_nf(x,u_n)dx
 =0.
 \end{eqnarray}
 Hence, by (\ref{3.1}), (\ref{3.3}), (\ref{3.8}), $(\phi_1)$, $(V)$ and (\ref{3.52}), we have
  \begin{eqnarray*}
  c& = & \lim_{n\rightarrow\infty}\left\{I(u_n)-\left\langle I'(u_n), \frac{1}{l}u_n\right\rangle\right\}\\
   & = & \lim_{n\rightarrow\infty}\left\{\iint_{\mathbb{R}^{2N}}\left(\Phi(|D_su_n|)
         -\frac{1}{l}a(|D_su_n|)|D_su_n|^2\right)d\mu\right.\\
   &   & +\int_{\mathbb{R}^N}V(x)\left(\Phi(|u_n|)-\frac{1}{l}a(|u_n|)u_n^2\right)dx\\
   &   & \left.+ \int_{\mathbb{R}^N}\left(\frac{1}{l}u_nf(x,u_n)-F(x,u_n)\right)dx\right\}\\
   &\leq& \lim_{n\rightarrow\infty}\left\{\int_{\mathbb{R}^N}\left(\frac{1}{l}u_nf(x,u_n)
         -F(x,u_n)\right)dx\right\}=0,
 \end{eqnarray*}
 which contradicts $c>0$. Therefore, $\lambda_3> 0$ and thus (\ref{3.48}) holds.
 \par
 Then, it follows from (\ref{3.48}) that there exist a constant $\delta>0$, a subsequence of $\{u_n\}$, still denoted by $\{u_n\}$, and a sequence $\{z_n\}\subset \mathbb{Z}^N$ such that
 \begin{eqnarray}\label{3.53}
 \int_{B_2(z_n)}\Phi(|u_n|)dx=\int_{B_2}(\Phi(|u_n^*|)dx>\delta, \quad \mbox{ for all } n\in \mathbb{N},
 \end{eqnarray}
 where $u_n^*:=u_n(\cdot +z_n)$. Since $V$ and $F$ are $1$-periodic in $x$, $\{u_n^*\}$ is also a $(C)_c$-sequence of $I$. Then, passing to a subsequence of $\{u_n^*\}$, still denoted by $\{u_n^*\}$, we can assume that there exists a $u^*\in W$ such that
 \begin{eqnarray}\label{3.54}
 u_n^*\rightharpoonup u^* \mbox{ in } W \quad \mbox{ and } \quad u_n^*\rightarrow u^* \mbox{ in } L^{\Phi}(B_2).
 \end{eqnarray}
 Thus, by (\ref{3.53}), (\ref{3.54}) and (\ref{ 2.2}), we obtain that $u^*\neq{\bf 0}$. Moreover, it follows from Lemma 3.7 and (\ref{3.8}) that
 \begin{eqnarray*}
  \|I'(u^*)\|_{W^*}\leq \liminf_{n\rightarrow\infty}\|I'(u_n^*)\|_{W^*}=0,
 \end{eqnarray*}
 which implies $I'(u^*)={\bf 0}$, that is, $u^*$ is a nontrivial solution of equation (\ref{eq1}). \quad $\Box$
 \vskip2mm
 \noindent
 {\bf Lemma 3.9.}\ \  {\it Assume that $(\phi_1)$, $(V)$ and $(f_1)$ hold. Then,}
  $$\langle I'(u),u\rangle=\langle I_1'(u),u\rangle-o(\langle I_1'(u),u\rangle)\quad\mbox{\it as }\quad\|u\|_{s,\Phi}\rightarrow 0.$$
 \vskip2mm
 \noindent
 {\bf Proof.}\ By using the continuity of $I_i'(i=1, 2)$ defined by (\ref{3.5+}) and (\ref{3.6+}), we can easily verify that $\langle I_i'(u),u\rangle=o(1)(i=1,2)$ as $\|u\|_{s,\Phi}\rightarrow 0$. Then, it is sufficient to prove $\langle I_2'(u),u\rangle=o(\langle I_1'(u),u\rangle)$ as $\|u\|_{s,\Phi}\rightarrow 0$ because $\langle I'(u),u\rangle=\langle I_1'(u),u\rangle-\langle I_2'(u),u\rangle$.
 \par
 For any given constant $\varepsilon>0$, it follows $(f_1)$, $(\phi_1)$ and (\ref{ 2.11}) that there exists a constant $C_\varepsilon>0$ such that
 \begin{eqnarray}\label{3.55}
 |tf(x,t)|\leq \varepsilon\Phi(|t|)+C_\varepsilon\Phi_*(|t|), \mbox{ for all } (x,t)\in\mathbb{R}^N\times\mathbb{R}.
 \end{eqnarray}
 Then, by (\ref{3.6+}) and (\ref{3.55}), we have
  \begin{eqnarray}\label{3.58}
 |\langle I_2'(u),u\rangle|&\leq &\int_{\mathbb{R}^N}|uf(x,u)|dx\nonumber\\
 &\leq & \varepsilon\int_{\mathbb{R}^N}\Phi(|u|)dx+C_{\varepsilon}\int_{\mathbb{R}^N}\Phi_{*}(|u|)dx.
 \end{eqnarray}
 Moreover, by (\ref{3.5+}), $(\phi_1)$ and $(V)$, we have
 \begin{eqnarray}\label{3.59}
  \langle I_1'(u),u\rangle
  &=&\iint_{\mathbb{R}^{2N}}a(|D_su|)|D_su|^2d\mu+\int_{\mathbb{R}^N}V(x)a(|u|)u^2dx\nonumber\\
  &\geq& l\iint_{\mathbb{R}^{2N}}\Phi(|D_su|)d\mu+\alpha_1l\int_{\mathbb{R}^N}\Phi(|u|)dx.
 \end{eqnarray}
  Then, (\ref{3.58}), (\ref{3.59}), Lemma 2.2, $3)$ in Lemma 2.5, $1)$ in Lemma 2.6 and the fact that $1<m<l^*$ imply that
 \begin{eqnarray*}
         \lim_{\|u\|_{s,\Phi}\rightarrow 0}\frac{|\langle I_2'(u),u\rangle|}{\langle I_1'(u),u\rangle}
 &\leq &  \lim_{\|u\|_{s,\Phi}\rightarrow 0}
 \frac{\varepsilon\int_{\mathbb{R}^N}\Phi(|u|)dx+C_{\varepsilon}\int_{\mathbb{R}^N}\Phi_{*}(|u|)dx}
 {l\iint_{\mathbb{R}^{2N}}\Phi(|D_su|)d\mu+\alpha_1l\int_{\mathbb{R}^N}\Phi(|u|)dx}\\
 &\leq & \frac{\varepsilon}{\alpha_1l}+\lim_{\|u\|_{s,\Phi}\rightarrow 0}
 \frac{C_{\varepsilon}\int_{\mathbb{R}^N}\Phi_{*}(|u|)dx}
 {\min\{1,\alpha_1\}l\left(\iint_{\mathbb{R}^{2N}}\Phi(|D_su|)d\mu+\int_{\mathbb{R}^N}\Phi(|u|)dx\right)}\\
  &\leq & \frac{\varepsilon}{\alpha_1l}+\lim_{\|u\|_{s,\Phi}\rightarrow 0}
 \frac{C_{\varepsilon}\max\{C_{\Phi_{*}}^{l^*},C_{\Phi_{*}}^{m^*}\}\|u\|_{s,\Phi}^{l^*}}
 {\min\{1,\alpha_1\}lC_m\|u\|_{s,\Phi}^{m}}\\
 &=& \frac{\varepsilon}{\alpha_1l}.
  \end{eqnarray*}
  Since $\varepsilon$ is arbitrary, we conclude that $|\langle I_2'(u),u\rangle|=o(\langle I_1'(u),u\rangle)$ as $\|u\|_{s,\Phi}\rightarrow 0$, which implies that $\langle I_2'(u),u\rangle=o(\langle I_1'(u),u\rangle)$ as $\|u\|_{s,\Phi}\rightarrow 0$. \quad $\Box$

  \vskip2mm
 \noindent
 {\bf Proof of Theorem 1.1 and Theorem 1.2} \quad Lemma 3.8 shows that equation (\ref{eq1}) has at least a nontrivial solution under the assumptions of Theorem 1.1 and Theorem 1.2, respectively. Next, we prove equation (\ref{eq1}) has a ground state solution. Let
 $$\mathcal{N}:=\{u\in W\setminus\{{\bf 0}\}: I'(u)={\bf 0}\} \quad \mbox{ and }\quad  d:=\inf_{u\in\mathcal{N}}\{I(u)\}.$$
 \par
 First, we claim that $d\geq0$. Indeed, for any given nontrivial critical point $u\in\mathcal{N}$, by (\ref{3.1}), (\ref{3.3}), $(\phi_1)$, $(V)$ and $(f_3)$ (or $(f_5)$), we have
 \begin{eqnarray*}
I(u)&=& I(u)-\left\langle I'(u), \frac{1}{m}u\right\rangle\nonumber\\
     &  = & \iint_{\mathbb{R}^{2N}}\left(\Phi(|D_su|)-\frac{1}{m}a(|D_su|)|D_su|^2\right)d\mu\nonumber\\
     &    & +\int_{\mathbb{R}^N}V(x)\left(\Phi(|u|)-\frac{1}{m}a(|u|)u^2\right)dx\nonumber\\
     &    & + \int_{\mathbb{R}^N}\left(\frac{1}{m}uf(x,u)-F(x,u)\right)dx\nonumber\\
     &\geq& \frac{1}{m}\int_{\mathbb{R}^N}\widehat{F}(x,u)dx\geq 0.
   \end{eqnarray*}
 Since the nontrivial critical point $u$ of $I$ is arbitrary, we conclude $d\geq 0$. Choose a sequence $\{u_n\}\subset\mathcal{N}$ such that $I(u_n)\rightarrow d$ as $n\rightarrow \infty$. Then, it is obvious that $\{u_n\}$ is a $(C)_{d}$-sequence of $I$ for the level $d$. Lemma 3.4 and Lemma 3.6 show that $\{u_n\}$ is bounded in $W$. Moreover, combining Lemma 3.9 with the fact that $\{u_n\}\subset\mathcal{N}$, we can conclude that there exists a constant $M_3>0$ such that
 \begin{eqnarray}\label{3.60}
 \|u_n\|_{s,\Phi}\geq M_3, \quad \mbox{ for all } n\in\mathbb{N}.
   \end{eqnarray}
 \par
 Now, we claim that
 \begin{eqnarray}\label{3.61}
 \lambda_4:=\lim_{n\rightarrow \infty}\sup_{y\in\mathbb{R}^N}\int_{B_2(y)}\Phi(|u_n|)dx>0.
 \end{eqnarray}
 Indeed, if $\lambda_4=0$, similar to (\ref{3.52}), we can get
  \begin{eqnarray}\label{3.62}
 \lim_{n\rightarrow\infty}\int_{\mathbb{R}^N}u_nf(x,u_n)dx=0.
 \end{eqnarray}
 Then, by (\ref{3.3}), $(\phi_1)$, $(V)$ and (\ref{3.62}), we have
  \begin{eqnarray*}
  0 &  = & \lim_{n\rightarrow\infty}\left\{\langle I'(u_n),u_n\rangle+\int_{\mathbb{R}^N}u_nf(x,u_n)dx\right\}\\
    &  = & \lim_{n\rightarrow\infty}\left\{\iint_{\mathbb{R}^{2N}}a(|D_su_n|)|D_su_n|^2d\mu
          +\int_{\mathbb{R}^N}V(x)a(|u_n|)u_n^2dx\right\}\\
   &\geq& \lim_{n\rightarrow\infty}\left\{l\iint_{\mathbb{R}^{2N}}\Phi(|D_su_n|)d\mu
          +\alpha_1l\int_{\mathbb{R}^N}\Phi(|u_n|)dx\right\}\\
   &\geq&0,
 \end{eqnarray*}
 which together with Lemma 2.2 implies that $\|u_n\|_{s,\Phi}=\|u_n\|_{\Phi}+[u_n]_{s, \Phi}\rightarrow 0$ as $n\rightarrow \infty$, which contradicts (\ref{3.60}). Therefore, $\lambda_4> 0$ and thus (\ref{3.61}) holds.
 \par
 Next, with similar arguments as those in Lemma 3.8, let $u_n^*:=u_n(\cdot +z_n)$. Then, $\{u_n^*\}$ is also a $(C)_d$-sequence of $I$. Moreover, there exist a subsequence of $\{u_n^*\}$, still denoted by $\{u_n^*\}$, and a $u^*\in W$ such that $u_n^*\rightharpoonup u^*$ in $W$ with $u^*\neq{\bf 0}$ and $I'(u^*)={\bf 0}$. This shows that $u^*\in\mathcal{N}$, and thus $I(u^*)\geq d$.
 \par
 On the other hand, by (\ref{3.1}), (\ref{3.3}), $(\phi_1)$, $(V)$, $(f_3)$ (or $(f_5)$) and Fatou's Lemma, we have
 \begin{eqnarray*}
I(u^*)&=& I(u^*)-\left\langle I'(u^*), \frac{1}{m}u^*\right\rangle\nonumber\\
     &  = & \iint_{\mathbb{R}^{2N}}\left(\Phi(|D_su^*|)-\frac{1}{m}a(|D_su^*|)|D_su^*|^2\right)d\mu\nonumber\\
     &    & +\int_{\mathbb{R}^N}V(x)\left(\Phi(|u^*|)-\frac{1}{m}a(|u^*|){|u^*|}^2\right)dx\nonumber\\
     &    & + \int_{\mathbb{R}^N}\left(\frac{1}{m}u^*f(x,u^*)-F(x,u^*)\right)dx\nonumber\\
     &\leq& \liminf_{n\rightarrow\infty}\left\{I(u_n^*)-\left\langle I'(u_n^*), \frac{1}{m}u_n^*\right\rangle\right\}\nonumber\\
     & = & d.
   \end{eqnarray*}
Therefore, $I(u^*)=d$, that is, $u^*$ is a ground state solution of equation (\ref{eq1}). This finishes the proof. \quad $\Box$

 \vskip-2mm
 \noindent
 \section{Examples}
 \setcounter{equation}{0}
 \vskip0mm
 \ \quad  For equation (\ref{eq1}), given $s\in(0, 1)$ and $N\in\mathbb{N}$, function $\phi$ defined by (\ref{1.2}) can be chosen from the following cases which satisfy conditions $(\phi_1)$-$(\phi_2)$:\\
 {\bf Case 1.} Let $\phi(t)=|t|^{p-2}t$ for $t\neq0$, $\phi(0)=0$ with $1<p<\frac{N}{s}$. In this case, simple computations show that $l=m=p$;\\
 {\bf Case 2.} Let $\phi(t)=|t|^{p-2}t+|t|^{q-2}t$ for $t\neq 0$, $\phi(0)=0$ with $1<p<q<\frac{N}{s}<\frac{pq}{q-p}$. In this case, simple computations show that $l=p, m=q$;\\
 {\bf Case 3.} Let $\phi(t)=\frac{|t|^{q-2}t}{\log(1+|t|^p)}$ for $t\neq 0$, $\phi(0)=0$ with $1<p+1<q<\frac{N}{s}<\frac{q(q-p)}{p}$. In this case, simple computations show that $l=q-p, m=q$.
 \par
 Moreover, we also give a case that satisfies condition $(\phi_1)$ but does not satisfy condition $(\phi_2)$:\\
 {\bf Case 4.} Let $\phi(t)=|t|^{q-2}t\log(1+|t|^p)$ for $t\neq 0$, $\phi(0)=0$ with $1<q<p+q<\frac{N}{s}<\frac{q(p+q)}{p}$. In this case, simple computations show that $l=q, m=p+q$.
 \par
 For example, regarding to Case 2, the operator in nonlocal problem (\ref{eq1}) defined by (\ref{1.3}) reduces to the following fractional $(p, q)$-Laplacian operator
 \begin{eqnarray*}
 (-\Delta_{p}-\Delta_{q})^{s}u(x)
 & = & P.V.\int_{\mathbb{R}^N}\frac{|u(x)-u(y)|^{p-2}(u(x)-u(y))}{|x-y|^{N+ps}}dy\\
 &&+P.V.\int_{\mathbb{R}^N}\frac{|u(x)-u(y)|^{q-2}(u(x)-u(y))}{|x-y|^{N+qs}}dy.
  \end{eqnarray*}
 \par
 Let $f(x,t)=qh(x)|t|^{q-2}t\log(1+|t|)+\frac{h(x)|t|^{q-1}t}{1+|t|}$, where $h\in C(\mathbb{R}^N, (0, +\infty))$ is 1-periodic in $x$. Then, $F(x,t)=h(x)|t|^{q}\log(1+|t|)$ and $\widehat{F}(x,t)=\frac{h(x)|t|^{q+1}}{1+|t|}$. It is easy to check that $f$ satisfies $(f_1)$-$(f_2)$, but does not satisfy the (AR) type condition $\mbox{(AR)}^*$. However, we can see that it satisfies $(f_3)$. Indeed, since $\frac{N}{s}<\frac{pq}{q-p}$, then there exists constant $k\in(\frac{N}{sp}, \frac{q}{q-p})$ such that
  \begin{eqnarray*}
 \limsup_{|t|\rightarrow \infty}\left(\frac{|F(x,t)|}{|t|^{l}}\right)^k\frac{1}{\widehat{F}(x,t)}
 =\limsup_{|t|\rightarrow \infty}\frac{h^{k-1}(x)(1+|t|)(\log(1+|t|))^k}{|t|^{(p-q)k+q+1}}=0,
  \end{eqnarray*}
 which implies that condition $(f_3)$ holds. Therefore, by using Theorem 1.1, we obtain that equation (\ref{eq1}) has at least one ground state solution when potential $V$ satisfies condition $(V)$.
 \par
 In addition, let $f(x,t)=h(x)\gamma(t)$, where $h\in C(\mathbb{R}^N, (0, +\infty))$ is 1-periodic in $x$ and
 \begin{eqnarray*}
 \gamma(t)=
 \begin{cases}
  \begin{array}{ll}
 0, &|t|\leq1,\\
 \left(|t|^{\frac{q+p^*-4}{2}}-\frac{1}{|t|}\right)t, &|t|>1.\\
    \end{array}
 \end{cases}
 \end{eqnarray*}
 Then, $F(x,t)=h(x)\Gamma(t)$, where
 \begin{eqnarray*}
 \Gamma(t)=
 \begin{cases}
  \begin{array}{ll}
 0, &|t|\leq1,\\
 \frac{2}{q+p^*}|t|^{\frac{q+p^*}{2}}-|t|+\frac{q+p^*-2}{q+p^*}, &|t|>1.\\
    \end{array}
 \end{cases}
 \end{eqnarray*}
 It is easy to check that $f$ satisfies $(f_1)$ and $(f_4)$, but does not satisfy $(f_3)$ and the (Ne) type condition $\mbox{(Ne)}^*$. However, we can see that it satisfies $(f_5)$. Indeed, since
 \begin{eqnarray}\label{4.1}
 \frac{1-\theta^l}{m}tf(x,t)=\frac{1-\theta^p}{q}h(x)t\gamma(t) \quad\mbox{and}\quad F(x,t)-F(x,\theta t)\leq F(x,t)=h(x)\Gamma(t),
   \end{eqnarray}
 for all $\theta\in\mathbb{R}, (x,t)\in\mathbb{R}^N\times\mathbb{R}.$ Then, it is obvious that
  \begin{eqnarray}\label{4.2}
\frac{1-\theta^l}{m}tf(x,t)\geq F(x,t)-F(x,\theta t), \mbox{ for all } \theta\in\mathbb{R}, (x,t)\in\mathbb{R}^N\times[-1,1].
   \end{eqnarray}
 Moreover,
$$\inf_{|t|>1}\frac{t\gamma(t)-q\Gamma(t)}{t\gamma(t)}= \inf_{|t|>1}\frac{\frac{p^*-q}{q+p^*}|t|^{\frac{q+p^*}{2}}+(q-1)|t|-\frac{q^2+qp^*-2q}{q+p^*}}
       {|t|^{\frac{q+p^*}{2}}-|t|}>0,$$
which implies that there exists a $\theta_0\in(0, 1)$ such that
  \begin{eqnarray}\label{4.3}
\frac{1-\theta^p}{q}h(x)t\gamma(t)\geq h(x)\Gamma(t), \mbox{ for all } \theta\in[0, \theta_0], x\in\mathbb{R}^N, |t|>1.
   \end{eqnarray}
Then, combining (\ref{4.2}) and (\ref{4.3}) with (\ref{4.1}), we can conclude that $(f_5)$ holds. Therefore, by using Theorem 1.2, we obtain that equation (\ref{eq1}) has at least one ground state solution when potential $V$ satisfies condition $(V)$.

 \section*{Funding information}
  This project is partially supported by the Guangdong Basic and Applied Basic Research Foundation (No: 2020A1515110706), Research Startup Funds of DGUT (No: GC300501-100), Yunnan Fundamental Research Projects (No: 202301AT070465), Xingdian Talent Support Program for Young Talents of Yunnan Province.

\end{document}